\documentclass{amsart}

\usepackage{amsthm}
\usepackage{amsfonts}
\usepackage{amsmath}
\usepackage{amssymb}
\usepackage{mathrsfs}
\usepackage{fullpage}
\usepackage{stmaryrd}
\usepackage{hyperref}
\usepackage{verbatim}

\theoremstyle{plain}

\theoremstyle{definition}

\theoremstyle{remark}

\newcommand{\Begin}[2]{\begin{#1}\label{#2}}

\newcommand{\Q}{\mathbb{Q}}

\newcommand{\bPi}{\mathbf{\Pi}}
\newcommand{\bSigma}{\mathbf{\Sigma}}
\newcommand{\bDelta}{\mathbf{\Delta}}

\newcommand{\bbP}{\mathbb{P}}
\newcommand{\bbQ}{\mathbb{Q}}

\newcommand{\CA}{\mathcal{A}}
\newcommand{\CB}{\mathcal{B}}
\newcommand{\CM}{\mathcal{M}}
\newcommand{\CN}{\mathcal{N}}
\newcommand{\CL}{\mathcal{L}}
\newcommand{\SCRL}{\mathscr{L}}
\newcommand{\SCRU}{\mathscr{U}}

\newcommand{\forces}{\Vdash}
\newcommand{\analytic}{{\bSigma_1^1}}
\newcommand{\lanalytic}{{\Sigma_1^1}}

\newcommand{\borel}{{\bDelta_1^1}}
\newcommand{\lborel}{{\Delta_1^1}}
\newcommand{\cantorspace}{{{}^\omega 2}}
\newcommand{\bairespace}{{{}^\omega\omega}}
\newcommand{\finBinarySequence}{{{}^{<\omega}2}}
\newcommand{\finNaturalSequence}{{{}^{<\omega}\omega}}

\newcommand{\F}{{F_{\omega_1}}}
\newcommand{\KP}{{\mathsf{KP}}}
\newcommand{\KPI}{{\mathsf{KP + INF}}}
\newcommand{\HYP}{{\mathcal{HYP}}}

\begin{document}

\title{The Countable Admissible Ordinal Equivalence Relation}

\author{William Chan}
\address{Department of Mathematics, California Institute of Technology, Pasadena, CA 91106}
\email{wcchan@caltech.edu}

\begin{abstract}
Let $\F$ be the countable admissible ordinal equivalence relation defined on $\cantorspace$ by  $x \ \F \ y$ if and only if $\omega_1^x = \omega_1^y$. Some invariant descriptive set theoretic properties of $F_{\omega_1}$ will be explored using infinitary logic in countable admissible fragments as the main tool. Marker showed $\F$ is not the orbit equivalence relation of a continuous action of a Polish group on $\cantorspace$. Becker stengthened this to show $\F$ is not even the orbit equivalence relation of a $\borel$ action of a Polish group. However, Montalban has shown that $\F$ is $\borel$ reducible to an orbit equivalence relation of a Polish group action, in fact, $\F$ is classifiable by countable structures. It will be shown here that $\F$ must be classified by structures of high Scott rank. Let $E_{\omega_1}$ denote the equivalence of order types of reals coding well-orderings. If $E$ and $F$ are two equivalence relations on Polish spaces $X$ and $Y$, respectively, $E \leq_{\text{a}\borel} F$ denotes the existence of a $\borel$ function $f : X \rightarrow Y$ which is a reduction of $E$ to $F$, except possibly on countable many classes of $E$. Using a result of Zapletal, the existence of a measurable cardinal implies $E_{\omega_1} \leq_{\text{a}\borel} \F$. However, it will be shown that in G\"{o}del's constructible universe $L$ (and set generic extensions of $L$), $E_{\omega_1} \leq_{\text{a}\borel} \F$ is false. Lastly, the techniques of the previous result will be used to show that in $L$ (and set generic extensions of $L$), the isomorphism relation induced by a counterexample to Vaught's conjecture cannot be $\borel$ reducible to $\F$. This shows the consistency of a negative answer to a question of Sy-David Friedman.
\end{abstract}

\maketitle\let\thefootnote\relax\footnote{January 28, 2016
\\*\indent Research partially supported by NSF grants DMS-1464475 and EMSW21-RTG DMS-1044448}


\section{Introduction}\label{Introduction}
If $x \in \cantorspace$, $\omega_1^x$ denotes the supremum of the order types of $x$-recursive well-orderings on $\omega$. Moreover, $\omega_1^x$ is also the minimum ordinal height of admissible sets containing $x$ as an element. The latter definition will be more relevant for this paper. 

The eponymous countable admissible ordinal equivalence relation, denoted by $\F$, is defined on $\cantorspace$ by:
$$x \ \F \ y \Leftrightarrow \omega_1^x = \omega_1^y$$

It is an $\lanalytic$ equivalence relation with all classes $\borel$. Moreover, $\F$ is a thin equivalance relation, i.e., it has no perfect set of inequivalence elements. Some further properties of $\F$ as an equivalence relation will be established in this paper.

Some basic results in admissibility theory and infinitary logic that will be useful throughout the paper will be reviewed in Section \ref{Admissibility and Infinitary Logic}. This section will cover briefly topics such as $\KP$, admissible sets, Scott ranks, and the Scott analysis. In this section, aspects of Barwise's theory of infinitary logic in countable admissible fragments, which will be the main tool in many arguments, will be reviewed. As a example of an application, a proof of a theorem of Sacks (Theorem \ref{Sacks theorem}), which establishes that every countable admissible ordinal is of the form $\omega_1^x$ for some $x \in \cantorspace$, will be given. This proof serves as a template for other arguments. Sacks theorem also explains why it is appropriate to call $\F$ the ``countable admissible ordinal equivalence relation". 

There have been some early work on whether $\F$ satisfies certain properties of equivalence relations related to generalization of Vaught's conjecture. For example, Marker in \cite{An-Analytic-Equivalence-Relation-Not-Arising} has shown that $\F$ is not induced by a continuous action of a Polish group on the Polish space $\cantorspace$. Becker in \cite{Topological-Vaught-Conjecture-and-Mininal-Counterexamples}, page 782, strenghened this to show that: the equivalence relation $\F$ is not an orbit equivalence relation of a $\borel$ group action of a Polish group. A natural question following these results would be whether $\F$ is $\borel$ reducible to equivalence relations induced by continuous or $\borel$ actions of Polish groups. If such reductions do exist, another question could be what properties must these reductions have.
\\*
\\*\indent In Section \ref{Classifiable by Countable Structure}, $\F$ will be shown to be $\borel$ reducible to a continuous action of $S_\infty$, i.e., it is classifiable by countable structures. An explicit $\lborel$ classification of $\F$ by countable structures in the language with a single binary relation symbol, due to Montalb\'{a}n, will be provided. The classification of $\F$ will use an effective construction of the Harrison linear ordering. This classification, denoted $f$, has the additional property that for all $x \in \cantorspace$, $\text{SR}(f(x)) = \omega_1^x + 1$. This example was provided by Montalb\'{a}n through communication with Marks and the author.

The explicit classification, $f$, mentioned above has images that are structures of high Scott rank. In Section \ref{Finer Aspects of Classification by Countable Structures}, it will be shown that this is a necessary feature of all classification of $\F$ by countable structures. The lightface version of the main result of this section is the following: 
\\*
\\*\noindent\textbf{Theorem \ref{classification F high scott rank}}
\textit{Let $\SCRL$ be a recursive language. Let $S(\SCRL)$ denote the set of reals that code $\SCRL$-structures on $\omega$. If $f : \cantorspace \rightarrow S(\SCRL)$ is a $\lborel$ function such that $x \ \F \ y$ if and only if $f(x) \cong_\SCRL f(y)$, then for all $x$, $\text{SR}(f(x)) \geq \omega_1^x$.}
\\*
\\*\noindent The more general form considers reductions that are $\lborel(z)$ and involves a condition on the admissible spectrum of $z$. Intuitively, Theorem \ref{classification F high scott rank} (in its lightface form as stated above) asserts that any potential classification of $\F$ must have high Scott rank in the sense that the image of any real under the reduction is a structure of high Scott rank. High Scott rank means that $\text{SR}(f(x))$ is either $\omega_1^x$ or $\omega_1^{x} + 1$.

Section \ref{Almost Borel Reductions} is concerned with a weak form of reduction of equivalence relations, invented by Zapletal, known as almost $\borel$ reduction. If $E$ and $F$ are two $\analytic$ equivalence relations on Polish space $X$ and $Y$, respectively, then $E$ is almost $\borel$ reducible to $F$ (in symbols: $E \leq_{\text{a}\borel} F$) if and only there is a $\borel$ function $f : X \rightarrow Y$ and a countable set $A$ such that if $x,y \notin A$, then $x \ E \ y$ if and only if $f(x) \ E \ f(y)$. 

An almost Borel reduction is simply a reduction that may fail on countably many classes. Often $\analytic$ equivalence relation may have a few unwieldly classes. The almost Borel reduction is especially useful since it can be used to ignore these classes. One example of such an $\analytic$ equivalence relation is $E_{\omega_1}$ which is the isomorphism relation of well-orderings with a single class of non-well-orderings. It is defined on $\cantorspace$ by:
\\*
\\* \text{ } \hfill $x \ E_{\omega_1} \ y \Leftrightarrow (x,y \notin WO) \vee (\text{ot}(x) = \text{ot}(y))$ \hfill \text{ }
\\*
\\*\noindent $E_{\omega_1}$ is a thin $\lanalytic$ equivalence with one $\Sigma_1^1$ class and all the other classes are $\borel$. 

Zapletal isolated an invariant of equivalence relations called the pinned cardinal. This invariant involves pinned names on forcings: an idea that appears implicitly or explicitly in the works of Silver, Burgess, Hjorth, and Zapletal in the study of thin $\analytic$ equivalence relations. Zapletal showed that there is a deep connection between $E_{\omega_1}$, almost $\borel$ reducibilities, and pinned cardinals under large cardinal assumptions:
\\*
\\*\noindent\textbf{Theorem \ref{Ew1 aleph0 almost reduction pinned cardinal}} (\cite{Forcing-Borel-Reducibily-Invariants} Theorem 4.1.3) 
\textit{If there exists a measurable cardinal and $E$ is a $\analytic$ equivalence relation with infinite pinned cardinal, then $E_{\omega_1} \leq_{\text{a}\borel} E$.}
\\*
\\*\indent Given that this result involves large cardinals, a natural question would be to explore the consistency results surrounding Zapletal's theorem. For example, a natural qustion is whether $\mathsf{ZFC}$ can prove the above result of Zapletal. More specifically, is this result true in G\"{o}del constructible universe $L$? This investigation leads to $\F$ in the following way: It will be shown that $\F$ has infinite pinned cardinal. Hence, with a measurable cardinal, $E_{\omega_1} \leq_{\text{a}\borel} \F$ via the result of Zapletal. (The author can show that $0^\sharp$ can prove the statement that $E_{\omega_1} \leq_{\text{a}\borel} \F$. A proof of this will appear in a future paper on pinned cardinals.) 

The main result of this section is
\\*
\\*\noindent\textbf{Theorem \ref{no almost reduction Ew1 F in L}} 
\textit{The statement $E_{\omega_1} \leq_{\text{a}\borel} \F$ is not true in $L$ (and set generic extensions of $L$).}
\\*
\\*\noindent This result is proved by using infinitary logic in admissible fragments to show that if $f$ is a $\lborel(z)$ function which witnesses $E_{\omega_1} \leq_{\text{a}\borel} \F$, then $z$ has an admissibility spectrum which is full of gaps relative to the set of all admissible ordinals. No constructible real (or even set generic over $L$ real) can have such property. 

The final section addresses a question of Sy-David Friedman using the techniques of the previous section. Essentially, the question is:
\\*
\\*\noindent \textbf{Question \ref{sy friedman vaught question}} Is it possible that the isomorphism relation of a counterexample to Vaught's conjecture is $\borel$ bireducible to $\F$?
\\*
\\*\indent The main result of this final section is:
\\*
\\*\noindent\textbf{Theorem \ref{counterexample vaught not reducible F}}
\textit{In $L$ (and set generic extensions of $L$), no isomorphism relation of a counterexample to Vaught's conjecture can be $\borel$ reducible to $\F$.}
\\*
\\*\indent This yields a negative answer to Friedman's question in $L$ and set generic extensions of $L$.

The author would like to acknowledge and thank Sy-David Friedman, Su Gao, Alexander Kechris, Andrew Marks, and Antonio Montalb\'{a}n for very helpful discussions and comments about what appears in this paper. 

\section{Admissibility and Infinitary Logic}\label{Admissibility and Infinitary Logic}

\indent The reader should refer to \cite{Admissible-Sets-and-Structures} for definitions and further details about admissibility. 

Let $\dot \in$ denote a binary relation symbol. Let $\SCRL$ be a language such that $\dot \in \in \SCRL$. $\KP_{\SCRL}$ denotes Kripke-Platek Set Theory in the language $\SCRL$ with $\dot\in$ serving as the distinguished membership symbol. The $\SCRL$ subscript will usually be concealed. $\KPI$ is $\KP$ augmented with the axiom of infinity.

\Begin{definition}{admissible sets}
Let $\SCRL$ be a language containing $\dot\in$. A $\SCRL$-structure $\CA = (A, \dot\in^\CA, ...)$ is an \textit{admissible set} if and only if $\CA \models \KP$, $A$ is a transitive set, and $\dot\in^\CA = \in \upharpoonright A$.

If $\CA$ is an admissible set, then $o(\CA) = A \cap \text{ON}$.

An ordinal $\alpha$ is an \textit{admissible ordinal} if and only if there is an admissible set $\CA$ such that $\alpha = o(\CA)$. More generally, if $x \in \cantorspace$, an ordinal $\alpha$ is $x$-admissible if and only if there is an admissible set $\CA$ such that $x \in A$ and $\alpha = o(\CA)$.

The \textit{admissibility spectrum} of $x$ is $\Lambda(x) = \{\alpha : \alpha \text{ is an $x$-admissible ordinal}\}$.

If $x \in \cantorspace$, $O(x) = \min(\Lambda(x))$.
\end{definition}

\Begin{definition}{hyperarithmetic structure}
For $x \in \cantorspace$, let $\mathcal{HYP}(x)$ denotes the $\subseteq$-smallest admissible set containing $x$.
\end{definition}

\Begin{definition}{w1x}
If $x \in \cantorspace$, let $\omega_1^x = O(x)$.
\end{definition}

\Begin{proposition}{sigma definability of L}
The function $(\alpha, x) \rightarrow L_{\alpha}(x)$, where $\alpha \in \text{ON}$ and $x$ is a set, is a  $\Sigma_1$ function in $\KP$. In fact, it is $\Delta_1$.
\end{proposition}

\begin{proof}
See \cite{Admissible-Sets-and-Structures}, Chapter II, Section 5 - 7. Also note that the function is defined on a $\Delta_1$ set. 
\end{proof}

\Begin{proposition}{L admissible ordinal is admissible}
If $\CA$ is an admissible set with $x \in A$ and $\alpha = o(\CA)$, then $L_\alpha(x)$ is an admissible set. In fact, $L_\alpha(x)$ is the $\subseteq$-smallest admissible set $\CA$ such that $x \in A$ and $o(\CA) = \alpha$.

In particular, if $\alpha$ is an $x$-admissible ordinal, then $L_\alpha(x)$ is an admissible set.
\end{proposition}

\begin{proof}
See \cite{Admissible-Sets-and-Structures}, Theorem II.5.7.
\end{proof}

\Begin{proposition}{HYP and L}
If $x \subseteq \omega$, then $\mathcal{HYP}(x) = L_{O(x)}(x) = L_{\omega_1^x}(x)$.
\end{proposition}

\begin{proof}
See \cite{Admissible-Sets-and-Structures}, Theorem II.5.9.
\end{proof}

\Begin{definition}{hyperarithmetic reals}
Let $x \in \cantorspace$. Suppose $\HYP(x) = (L_{\omega_1^x}(x), \in)$. Let $\text{HYP}^x = \cantorspace \cap L_{\omega_1^x}(x)$. $\text{HYP}^x$ is the set of all $x$-hyperarithmetic reals. 

In particular, $x$-hyperarithmetic reals are exactly those reals that appear in all admissible sets containing $x$.
\end{definition}

Next, the relevant aspects of first order infinitary logic and admissible fragments will be reviewed. The detailed formalization can be found in \cite{Admissible-Sets-and-Structures}, Chapter III.

\Begin{definition}{Lww and Linfw}
Let $\SCRL$ denote a first order language (a set of constant, relation, and function symbols). Fix a $\Delta_1$ class $\{v_\alpha : \alpha \in \text{ON}\}$, which will represent variables. $\SCRL_{\omega\omega}$ denotes the collection of finitary $\SCRL$-formulas using variables from $\{v_i : i < \omega\}$. $\SCRL_{\infty\omega}$ denotes the collection of all infinitary formulas with finitely many free variables. 
\end{definition}

\Begin{proposition}{definability Lww and Linfw}
In $\KPI$, $\SCRL_{\omega\omega}$ is a set. In $\KP$, $\SCRL_{\infty\omega}$ is a $\Delta_1$ class.
\end{proposition}

\begin{proof}
See \cite{Admissible-Sets-and-Structures}, Proposition III.1.4 and page 81.
\end{proof}

\Begin{proposition}{definability satisfaction}
($\KP$) ``$\CM \models_{\SCRL} \varphi(\bar{x})$'' as a relation on the language $\SCRL$, $\SCRL$-structure $\CM$, infinitary $\SCRL$-formula $\varphi$, and tuple $\bar{x}$ of $M$ is equivalent to a $\Delta_1$ predicate.
\end{proposition}

\begin{proof}
See \cite{Admissible-Sets-and-Structures}, pages 82-82.
\end{proof}

\Begin{definition}{admissible fragments}
Let $\SCRL$ be a language. Let $\CA$ be an admissible set such that $\SCRL$ is $\Delta_1$ definable in $\CA$. The admissible fragment of $\SCRL_{\infty\omega}$ given by $\CA$, denoted $\SCRL_\CA$, is defined as
$$\SCRL_\CA = \{\varphi \in A : \varphi \in L_{\infty\omega}\} = \{\varphi \in A : \CA \models \varphi \in L_{\infty\omega}\}$$
The last equivalence follows from $\Delta_1$ absoluteness.
\end{definition}

\Begin{definition}{well founded part and solidness}
Let $\SCRL$ be a language consisting of a binary relation $\dot\in$. Let $\CM$ be a $\SCRL$-structure such that $(M, \dot\in^\CM)$ satisfies extensionality. Define $\text{WF}(\CM)$ as the substructure consisting of the well-founded elements of $M$. $\text{WF}(\CM)$ is called the well-founded part of $\CM$. 

$\CM$ is called \textit{solid} if and only if $\text{WF}(\CM)$ is transitive.
\end{definition}

\Begin{remark}{solid remark}
The notion of solid comes from Jensen's \cite{Admissible-Sets}. Every structure has an isomorphic solid model that is obtained by Mostowski collapsing the well-founded part. 

The notion of solidness is mostly a convenience: In our usage, $\omega \subseteq M$. Therefore, Mostowski collapsing will not change reals. Transitivity is desired due to the definition of admissibility and in order to apply familiar absoluteness results. Rather than having to repeatedly Mostowski collapse $\text{WF}(\CM)$ and mention reals are not moved, one will just assume the well-founded part is transitive by demanding $\CM$ is solid. 
\end{remark}

\Begin{lemma}{truncation lemma}
(Truncation Lemma) If $\CM \models \mathsf{KP}$, then $\text{WF}(\CM) \models \KP$. In particular, if $\CM$ is a solid model, then $\text{WF}(\CM)$ is an admissible set.
\end{lemma}

\begin{proof}
See \cite{Admissible-Sets-and-Structures}, II.8.4.
\end{proof}

The following is the central technique used in the paper:

\Begin{theorem}{solid model existence theorem}
(Solid Model Existence Theorem) Let $\CA$ be a countable admissible set. Let $\SCRL$ be a language which is $\Delta_1$ definable over $\CA$ and contains a binary relation symbol $\dot\in$ and constant symbols $\bar{a}$ for each $a \in A$. Let $T$ be a consistent $\SCRL$-theory in the countable admissible fragment $\SCRL_\CA$, be $\Sigma_1$ definable over $\CA$, and contains the following:

(i) $\KP$

(ii) For each $a \in A$, the sentence $(\forall v)(v \dot\in \bar{a} \Rightarrow \bigvee_{z \in a} v =\bar{z})$. 

\noindent Then there exists a solid $\SCRL$-structure $\mathcal{B}$ such that $\mathcal{B} \models T$ and $\text{ON} \cap B = \text{ON} \cap A$.
\end{theorem}

\begin{proof}
See \cite{Admissible-Sets}, Section 4, Lemma 11.
\end{proof}

\Begin{theorem}{Sacks theorem}
(Sacks' Theorem) If $\alpha > \omega$ is an admissible ordinal, then there exists some $x \in \cantorspace$ such that $\alpha = \omega_1^x$.

Let $z \in \cantorspace$. If $\alpha \in \Lambda(z) \cap \omega_1$, then there exists $y \in \cantorspace$ with $\omega_1^y = \alpha$ and $z \leq_\text{T} y$.
\end{theorem}

\begin{proof}
See \cite{Countable-Admissible-Ordinals-and-Hyperdegrees}, Corollary 3.16. The following proof is similar to \cite{Admissible-Sets}, Section 4, Lemma 10. The second statement will be proved below:

Since $\alpha \in \Lambda(z)$, let $\CA$ be an admissible set such that $z \in A$ and $o(\CA) = \alpha$. (For example, $\CA = L_{\alpha}(z)$ by Proposition \ref{L admissible ordinal is admissible}.)

Let $\SCRL$ be a language consisting of the following:

(I) A binary relation symbol $\dot\in$.

(II) Constant symbols $\bar{a}$ for each $a \in A$. 

(III) One other distinguished constant symbol $\dot c$. 

\noindent The elements of $\SCRL$ can be appropriately coded as elements of $A$ so that $\SCRL$ is $\Delta_1$ definable over $\CA$. 

Let $T$ be a theory in the countable admissible fragment $\SCRL_{\CA}$ consisting of the following:

(i) $\KP$

(ii) For each $a \in A$, $(\forall v)(v \dot\in\bar{a} \Rightarrow \bigvee_{z \in a} v = \bar{z})$. 

(iii) $\dot c \subseteq \bar{\omega}$. 

(iv) For each ordinal $\sigma \in \alpha$, ``$\bar{\sigma}$ is not admissible relative to $\dot c$''. More formally, ``$L_{\bar\sigma}(\dot c) \not\models \KPI$''.

(v) $\bar{z} \leq_\text{T} \dot c$.

\noindent $T$ can be coded as a class in $A$ in such a way that it is $\Sigma_1$ in $\CA$. $T$ is consistent: Find any $u \in \cantorspace$ which codes an ordinal greater than $\alpha$. Let $c = u \oplus z$. Consider the following $\SCRL$-structure $\CM$: The universe $M$ is $H_{\aleph_1}$. For each $a \in \CA$, $\bar{a}^\CM = a$. (Since $A$ is countable and transitive, $A \in H_{\aleph_1}$.) $\dot\in^\CM = \in \upharpoonright H_{\aleph_1}$. $\dot c^\CM = c$. $\CM$ clearly satisfy (i), (ii), (iii), and (v). For (iv), suppose there is an $\sigma < \alpha$ such that $L_\sigma(c) \models \KP$. Since $c \in L_\sigma(c)$ and $L_\sigma(c) \models \KP$, $u \in L_\sigma(c)$ because $c = u \oplus z$. Since the Mostowski collapse map is a $\Sigma_1$ definable function in $\KP$, if reals code binary relations in the usual way, then $\KP$ proves the existence of $\text{ot}(u)$. Thus $\text{ot}(u) \in L_\sigma(c)$. However, $\text{ot}(u) > \alpha > \sigma$. Contradiction. It has been shown that $\CM$ also satisfy (iv). $T$ is consistent.

The Solid Model Existence Theorem (Theorem \ref{solid model existence theorem}) implies there is a solid $\SCRL$-structure $\mathcal{B} \models T$ such that $\text{ON} \cap B = \text{ON} \cap A = \alpha$. Let $y = \dot c^\CB$. The claim is that $\omega_1^y = \alpha$. By Lemma \ref{truncation lemma}, $\text{WF}(\CB)$ is an admissible set containing $y$ and $z$. $o(\text{WF}(\CB)) = \text{ON} \cap \text{WF}(\CB) = \text{ON} \cap B = \text{ON} \cap A = \alpha$. Thus $\omega_1^y \leq \alpha$. Now suppose that $\omega_1^y < \alpha$. In $V$, $L_{\omega_1^y}(y) \models \KP$. Since the function $(\alpha, x) \mapsto L_\alpha(x)$ is $\Delta_1$ (by Proposition \ref{sigma definability of L}) and the satisfaction relation is $\Delta_1$ (by Proposition \ref{definability satisfaction}), by $\Delta_1$ absoluteness between the transitive sets $\text{WF}(\CB)$ and $V$, one has $\text{WF}(\CB) \models L_{\omega_1^y}(y) \models \mathsf{KP}$. Again by absoluteness of $\Delta_1$ formulas between the transitive (in the sense of $\CB$) sets $\text{WF}(\CB)$ and $\CB$, $\CB \models L_{\omega_1^x} \models \KP$. Letting $\sigma = \omega_1^x < \alpha$, $\CB \models L_{\bar\sigma}(\dot c) \models \KP$. This contradicts $\CB \models T$. A similar absoluteness argument shows that $z \leq_\text{T} y$.
\end{proof}

\Begin{remark}{remark Sacks theorem}
This proof of Sacks theorem is the basic template for several other arguments throughout the paper. This proof will be frequently referred.
\end{remark}

Next, various aspects of the Scott analysis will be reviewed. Since there are some minor variations among the definitions of Scott rank, Scott sentences, canonical Scott sentences, etc., these will be provided below. See \cite{Model-Theory-an-Introduction}, page 57-60 or \cite{Scott-Sentences-and-Admissible-Sets} for more information.

\Begin{definition}{scott relation and sentence}
Let $\SCRL$ be a language. Define the binary relation $(\CM, \underline{a}) \sim_\alpha (\CN, \underline{b})$ where $\alpha \in \text{ON}$, $\underline{a} \in {}^{<\omega}M$, and $\underline{b} \in {}^{<\omega} N$ as follows:

(i) $(\CM, \underline{a}) \sim_0 (\CN, \underline{b})$ if and only if for all atomic $\SCRL$-formulas $\varphi$, $\CM \models \varphi(\underline{a})$ if and only if $\CN \models \varphi(\underline{b})$.

(ii) If $\alpha$ is a limit ordinal, then $(\CM, \underline{a}) \sim_\alpha (\CN, \underline{b})$ if and only if for all $\beta < \alpha$, $(\CM, \underline{a}) \sim_\beta (\CN, \underline{b})$.

(iii) If $\alpha = \beta + 1$, then $(\CM, \underline{a}) \sim_\alpha (\CN, \underline{b})$ if and only if for all $c \in M$, there exists a $d \in N$ such that $(\CM, \underline{a},c) \sim_\beta (\CN, \underline{b}, d)$ and for all $d \in N$, there exists a $c \in M$ such that $(\CM, \underline{a}, c) \sim_\beta (\CN, \underline{b},d)$.
\\*
\\*\indent Let $\CM$ be a $\SCRL$-structure and $\underline{a} \in {}^kM$ for some $k \in \omega$. For $\alpha \in \text{ON}$, the $\SCRL_{\infty\omega}$ formula $\Phi_{\underline{a},\alpha}^\CM(\underline{v})$ (in variables $\underline{v}$ such that $|\underline{v}| = k$) is defined as follows:

(I) Let $X$ be the set of all atomic and negation atomic $\SCRL$-formulas with free variables $\underline{v}$ such that $|\underline{v}| = k$. Let $\Phi_{\underline{a},0}^\CM (\underline{v}) = \bigwedge X$

(II) If $\alpha$ is a limit ordinal, let $X = \{\Phi_{\underline{a}, \beta}^\CM(\underline{v}) : \beta < \alpha\}$. Let $\Phi_{\underline{a},\alpha}^\CM(\underline{v}) = \bigwedge X$.

(III) If $\alpha = \beta + 1$, then let $X = \{(\exists w)\Phi_{\underline{a}b,\beta}^\CM(\underline{v}, w) : b \in M\}$ and $Y = \{\Phi_{\underline{a}b,\beta}^\CM(\underline{v},w): b \in M\}$. Then let $\Phi_{\underline{a},\alpha}^\CM(\underline{v}) = \bigwedge X \wedge (\forall w)\bigvee Y$.
\\*
\\*\indent For $\CM$, a $\CL$-structure, and $\underline{a} \in {}^kM$ (for some $k$), define $\rho(\CM,\bar{a})$ to be the least $\alpha \in \text{ON}$ such that for all $\underline{b} \in {}^kM$, $(\CM, \underline{a}) \sim_\alpha (\CM, \underline{b})$ if and and only if for all $\beta$, $(\CM, \underline{a}) \sim_\beta (\CM, \underline{b})$. 

Define $\text{SR}(\CM) = \sup\{\rho(\CM, \underline{a}) + 1 : \underline{a} \in {}^{<\omega}M\}$. Define $\text{R}(\CM) = \sup\{\rho(\CM, \underline{a}) : \underline{a} \in {}^{<\omega}M\}$.

Let $\alpha = \text{R}(\CM)$. Let 
$$X = \{(\forall \bar{v})(\Phi_{\underline{a},\alpha}^\CM(\underline{v}) \Rightarrow \Phi_{\underline{a},\alpha+1}^\CM(\underline{v})) : \underline{a} \in {}^{<\omega}M\}$$
$$\text{CSS}(\CM) = \Phi_{\emptyset, \alpha}^\CM \wedge \bigwedge X$$
$\text{CSS}(\CM)$ is the \textit{canonical Scott sentence} of $\CM$. $\text{SR}(\CM)$ is the \textit{Scott rank} of $\CM$.
\end{definition}

The following are well-known results. Usually, a careful inspection of the proof indicates what can be done in $\KPI$ or $\mathsf{ZFC}$.

\Begin{proposition}{scott relation KP}
The relation $\sim$ is equivalent to $\Delta_1$ formula over $\mathsf{KP + INF}$. 
\end{proposition}

\begin{proof}
It can be defined by $\Sigma$-recursion.
\end{proof}

\Begin{proposition}{phi formula}
Let $\mathcal{A}$ be an admissible set such that $\mathcal{A} \models \mathsf{INF}$. Let $\mathscr{L} \in A$ be a language. Let $\CM \in A$ range over $\SCRL$-structure, $\underline{a}$ range over elements of ${}^{<\omega}M$, and $\alpha$ range over $\text{ON} \cap A$. Then the function $f(\CM, \underline{a}, \alpha) = \Phi_{\underline{a}, \alpha}^\CM(\underline{v})$ is $\Delta_1$ definable in $\mathcal{A}$.

In particular, if $\text{R}(\CM) \in A$, then $\text{CSS}(\CM) \in A$.
\end{proposition}

\begin{proof}
It can be defined by $\Sigma$-recursion.
\end{proof}

\Begin{proposition}{sim and phi formula}
($\mathsf{KP + INF}$) Let $\SCRL$ be a language. Let $\CM$ and $\CN$ be $\SCRL$-structures, $\underline{a} \in {}^kM$, $\underline{b} \in {}^kN$ for some $k \in \omega$, and $\alpha \in \text{ON}$. Then $(\CM, \underline{a}) \sim_{\alpha} (\CN, \underline{b})$ if and only if $\CN \models \Phi_{\underline{a}, \alpha}^\CM(\underline{b})$.
\end{proposition}

\begin{proof}
This is proved by induction. See \cite{Model-Theory-an-Introduction}, Lemma 2.4.13.
\end{proof}

\Begin{proposition}{sim and elementary equivalence}
($\mathsf{KP + INF}$) If $(M, \underline{a}) \equiv_{\SCRL_{\infty\omega}} (N,\underline{b})$, then for all $\alpha$, $(M,\underline{a}) \sim_{\alpha} (N, \underline{b})$.
\end{proposition}

\begin{proof}
$\CM \models \varphi_{\underline{a},\alpha}(\underline{a})$. So $\CN \models \varphi_{\underline{a},\alpha}(\underline{b})$. By Proposition \ref{sim and phi formula}, $(\CM, \underline{a}) \sim_{\alpha} (\CN, \underline{b})$. 
\end{proof}

\Begin{definition}{quantifier rank}
Let $\SCRL$ be a language. Let $\varphi$ be a formula of $L_{\infty\omega}$. The \textit{quantifier rank} of $\varphi$ denoted $\text{qr}(\varphi)$ is defined as follows: 

(i) $\text{qr}(\varphi) = 0$ if $\varphi$ is an atomic formula. 

(ii) $\text{qr}(\neg \varphi) = \text{qr}(\varphi)$

(iii) $\text{qr}(\bigwedge X) = \text{qr}(\bigvee X) = \sup\{\text{qr}(\psi) : \psi \in X\}$.

(iv) $\text{qr}(\exists v\varphi) = \text{qr}(\forall v \varphi) =  \text{qr}(\varphi) + 1$. 
\end{definition}

\Begin{proposition}{qr definability}
The relation ``$\text{qr}(\varphi) = \alpha$'' is $\Delta_1$ definable in $\KPI$.
\end{proposition}

\begin{proof}
It can be defined by $\Sigma$-recursion.
\end{proof}

\Begin{proposition}{sim and quantifier rank}
($\mathsf{KP + INF}$) Let $\SCRL$ be a language. $\CM, \CN$ be $\SCRL$-structures. $\underline{a} \in {}^kM$ and $\underline{b} \in {}^kN$ for some $k \in \omega$. Then for all $\alpha \in \text{ON}$, $(\CM, \underline{a}) \sim_{\alpha} (\CN, \underline{b})$ if and only if for all $\varphi$ with $\text{qr}(\varphi) \leq \alpha$, $\CM \models \varphi(\underline{a})$ if and only if $\CN \models \varphi(\underline{b})$.
\end{proposition}

\begin{proof}
This is proved by induction. 
\end{proof}

\Begin{proposition}{elem equiv in admissible iff elem equiv}
($\mathsf{ZF}$) Let $\SCRL$ be some language. Let $\CM$ and $\CN$ be $\SCRL$-structures. Suppose $\CA$ is an admissible set with $\SCRL, \CM, \CN \in A$. Then $\CA \models \CM \equiv_{\SCRL_{\infty\omega}} \CN$ if and only if $\CM \equiv_{\SCRL_{\infty\omega}} \CN$.
\end{proposition}

\begin{proof}
See \cite{Scott-Sentences-and-Admissible-Sets}, Theorem 1.3.
\end{proof}

\Begin{remark}{elem equiv in admissible remark}
A common phenomenon is that certain properties are reflected between appropriate admissible sets and the true universe. A useful observation is that if such a property holds from the point of view of an admissible set then it is true in the universe. The above proposition asserts that infinitary elementary equivalence is such a property. 

Another familiar example is the effective boundedness theorem. Suppose $\varphi : \text{WO} \rightarrow \omega_1$ is a $\Pi_1^1$ rank. Let $B \subseteq \text{WO}$ be $\bSigma_1^1$. Let $\CA$ be a countable admissible set containing the parameters used to define $B$. Inside of $\CA$, $\varphi(B)$ is bounded by $o(\CA)$. A priori, the true bound on $\varphi(B)$ may be higher as the true universe has more countable ordinals and more members of $B$. However, the effective boundedness theorem asserts that in fact, in the true universe, $\varphi(B)$ is bounded by $o(\CA)$. 

The following proposition with an included proof shows countable admissible sets can also be used to produce true bounds on the Scott rank.
\end{remark}

\Begin{proposition}{length of scott rank}
Let $\SCRL$ be a countable language and $\CM$ be a countable $\SCRL$-structure. One may identify $\CM$ as a real by associating it with an isomorphic structure on $\omega$. If $\CA$ is an admissible set with $\SCRL, \CM \in A$, then $R(\CM) \leq \text{ON} \cap A$. $R(\CM) \leq O(\CM)$. $\text{SR}(\CM) \leq O(\CM) + 1$. In particular, $\text{R}(\CM) \leq \omega_1^\CM$ and $\text{SR}(\CM) \leq \omega_1^\CM + 1$. 
\end{proposition}

\begin{proof}
See \cite{Scott-Sentences-and-Admissible-Sets}, Corollary 1. 

It suffices to show that $R(M) \leq O(\CM)$. Suppose not. Then there exists $\underline{a}$ and $\underline{b}$ such that for all $\alpha < O(\CM)$, $(\CM, \underline{a}) \sim_\alpha (\CM, \underline{b})$ but for $\beta = O(\CM )$, $(\CM, \underline{a}) \not\sim_\beta (\CM, \underline{b})$. Let $\CA$ be an admissible set with $\CM \in A$ and $o(\CA) = O(\CM)$. By $\Delta_1$-absoluteness and Proposition \ref{sim and quantifier rank}, $A \models (\CM, \underline{a}) \equiv_{\SCRL_{\infty\omega}} (\CM, \underline{b})$. Thus by Proposition \ref{elem equiv in admissible iff elem equiv}, $(\CM, \underline{a}) \equiv_{\SCRL_{\infty\omega}} (\CM, \underline{b})$. However, $(\CM, \underline{a}) \not\sim_\beta (\CM, \underline{b})$ implies $\CM \models \Phi_{\underline{a}, \beta}^\CM(\underline{a})$ and $\CM \not\models \Phi_{\underline{a},\beta}^\CM(\underline{b})$ by Proposition \ref{sim and phi formula}. This shows $(\CM, \underline{a}) \not\equiv_{\SCRL_{\infty\omega}} (\CM, \underline{b})$. Contradiction.

\end{proof}

\Begin{definition}{scott sentence}
Let $\SCRL$ be a language. Let $\CM$ be a $\SCRL$-structure. $\varphi$ is a Scott sentence if and only if for all $\SCRL$-structure $\CN$ and $\CM$, $\CN \models \varphi$ and $\CM \models \varphi$ implies $\CM \equiv_{\SCRL_{\infty\omega}} \CN$.
\end{definition}

\Begin{theorem}{scott theorem}
($\mathsf{ZFC}$) Let $\SCRL$ be a language. Let $\CM$ be a countable $\SCRL$-structure. Then there exists a $L_{\infty\omega}$-sentence $\varphi$ such that for all countable $\SCRL$-structure $\CN$, $\CN \cong_{\SCRL} \CM$ if and only if $\CN \models \varphi$. In fact, $\varphi$ is $\text{CSS}(\CM)$. 

($\mathsf{KP + INF}$) If $\varphi$ is a Scott sentence for a countable structure $\CM$, then for all countable $\CN$, $\CN \models \varphi$ if and only if $\CN \cong_{\SCRL} \CM$.
\end{theorem}

\begin{proof}
Observe the first statement asserts that there exists a sentence such that whenever a countable structure satisfies this sentence, there exists an isomorphism between it and $\CM$. The existence of this sentence requires working beyond $\KPI$. The second statement asserts that $\KPI$ can prove that if a Scott sentence happens to exist, then for any countable structure satisfying this sentence, there is an isomorphism between it and $\CM$.

This is the Scott's isomorphism theorem. See \cite{Model-Theory-an-Introduction}, Theorem 2.4.15 for a proof. The results in $\KPI$ follows essentially the same proof with the assistance of some of the above propositions proved in $\KPI$.
\end{proof}

\Begin{definition}{space of structures}
Let $\SCRL$ be a countable language. Let $S(\SCRL)$ denote the set of all $\SCRL$-structures on $\omega$.
\end{definition}

\Begin{definition}{w-model KP}
Let $\dot \in$ be a binary relation symbol. Let $S^*$ denote the subset of $S(\{\dot\in\})$ consisting of $\omega$-models of $\KPI$.
\end{definition}

\Begin{proposition}{first order satisfaction borel}
Let $\{\phi_e : e \in \omega\}$ be a recursive enumeration of $\{\dot\in\}_{\omega\omega}$-formulas. The relation on $x \in S(\{\dot\in\})$ and $e \in \omega$ asserting ``$x \models \phi_e$'' is $\Delta_1^1$.

Also $S^*$ is $\Delta_1^1$.
\end{proposition}

\begin{proof}
See \cite{Recursive-Aspects-of-Descriptive-Set-Theory}, page 14-16 for relevant definitions and proofs.
\end{proof}

\Begin{remark}{reals in omega models}
One can check that there is a $\Delta_1^1$ function such that given $A \in S^*$ and $n \in \omega$, the function gives the element of $A$ which $A$ thinks is $n$. Using this, one can determine in a $\Delta_1^1$ way whether $A \in S^*$ thinks some $x \in \cantorspace$ exists. In the following, if $A \in S^*$ and $x \in \cantorspace$, the sentence ``$x \in A$'' should be understood as this informally described $\Delta_1^1$ relation.
\end{remark}

\Begin{proposition}{complexity of models of formula}
Let $\SCRL$ be a recursive language. Let $\varphi \in \HYP(x) \cap \SCRL_{\infty\omega}$. Then $\text{Mod}(\varphi) = \{s \in S(\SCRL) : s \models_{\SCRL}\varphi\}$ is $\Delta_1^1(x)$.
\end{proposition}

\begin{proof}
Note that $s \in \text{Mod}(\varphi)$ if and only if
$$(\exists A)(A \in S^* \wedge x \in A \wedge s \in A \wedge A \models s \models_{\SCRL}\varphi)$$
if and only if
$$(\forall A)(A \in S^* \wedge x \in A \wedge s \in A) \Rightarrow A \models s \models_{\SCRL} \varphi)$$
These equivalences are established using the absoluteness of satisfaction. This shows that $\text{Mod}(\varphi)$ is $\Delta_1^1(x)$.
\end{proof}

\Begin{remark}{philosophical remark on classification result}
Later, the paper will be concerned with relating countable admissible sets and isomorphism of countable structures. The second statement of Theorem \ref{scott theorem} captures the essence of these types of arguments: Isomorphism of countable structures is reflected between the true universe and admissible sets which witness the countability of the relevant structures and possesses a Scott sentence for these structures.

The original arguments for some results of this paper used more directly the second statement of Theorem \ref{scott theorem}. The argument presented below is simpler using the Scott isomorphism theorem and Proposition \ref{complexity of models of formula} but may conceal this essential idea.
\end{remark}

Now to introduce the main equivalence relation of this paper:

\Begin{definition}{F equivalence relation}
Let $\F$ be the equivalence relation defined on $\cantorspace$ by $x \ \F \ y$ if and only if $\omega_1^x = \omega_1^y$. $\F$ is a $\Sigma_1^1$ equivalence relation with all classes $\bDelta_1^1$.
\end{definition}

The first claim from the above definition is well known and follows easily from the characterization of $\omega_1^x$ as the supremum of the $x$-recursive ordinals. The next proposition implies each class is $\bDelta_1^1$. There will be much to say later about the complexity of each $\F$-equivalence class. 

\Begin{proposition}{F classes borel in real of higher church kleene}
Let $\alpha$ be a countable admissible ordinal and $z \in \cantorspace$ be such that $\alpha < \omega_1^z$. Then the set $\{y \in \cantorspace : \omega_1^y = \alpha\}$ is $\Delta_1^1(z)$.
\end{proposition}

\begin{proof}
If $u$ and $v$ are reals coding linear orderings on $\omega$, then $u \preceq v$ means there exists an order preserving injective function $f$ from the linear ordering coded by $u$ to the linear ordering coded by $v$. $\preceq$ is a $\Sigma_1^1$ relation in the variables $u$ and $v$.

Since $\alpha < \omega_1^z$, there exists some $e \in \omega$ such that $\{e\}^z$ is the characteristic function of a well-ordering isomorphic to $\alpha$. Let $B = \{y \in \cantorspace : \alpha = \omega_1^y\}$. Then
$$y \in B \Leftrightarrow (\forall n)\Big((\{n\}^y \in \text{WO} \Rightarrow \{n\}^y \preceq \{e\}^z )\Big) \wedge$$
$$(\forall k)(\exists j)(\{j\}^y \preceq \{e\}^z \wedge \{e\}^z\upharpoonright k \preceq \{j\}^y)$$
$B$ is $\Sigma_1^1(z)$. Also
$$y \notin B \Leftrightarrow (\exists j)(\forall n)(\{n\}^y \in \text{WO} \Rightarrow (\{n\}^y \preceq \{e\}^z \upharpoonright j) \vee (\exists n)(\{n\}^y \in \text{LO} \wedge \{n\}^y \preceq \{e\}^z \wedge \{e\}^z \preceq \{n\}^y)$$
$B$ is $\Pi_1^1(z)$. Hence $B$ is $\Delta_1^1(z)$.
\end{proof}

\section{Classifiable by Countable Structures}\label{Classifiable by Countable Structure}
\Begin{definition}{recursive pseudo-wellorderings}
Let $x \in \cantorspace$. A linear ordering $R$ on $\omega$ is an $x$-recursive $x$-\textit{pseudo-wellordering} if and only if $R$ is an $x$-recursive linear ordering on $\omega$ which is not a wellordering but $L_{\omega_1^x}(x) \models R$ is a wellordering, i.e. $R$ has no $x$-hyperarithmetic descending sequences. 
\end{definition}

\Begin{proposition}{existence pseudo-wellordering}
(Harrison, Kleene) For all $x \in \cantorspace$, there exists an $x$-recursive $x$-pseudo-wellordering.
\end{proposition}

\begin{proof}
See \cite{On-the-Forms-Predicate-II} or \cite{Higher-Recursion-Theory}, III.2.1. A generalized form of this construction will be used below. 

This can also be proved using Theorem \ref{solid model existence theorem} and infinitary logic in admissible fragments. In the application of Theorem \ref{solid model existence theorem}, Barwise compactness is used to show the consistency of the appropriate theory in the countable admissible fragment. See Nadel's proof given in \cite{Model-Theoretic-Logics} VIII, Section 5.7 for more details.
\end{proof}

The following characterizes the order type of $x$-recursive $x$-pseudo-wellorderings:

\Begin{theorem}{harrison theorem}
(Harrison) Let $R$ be a $x$-recursive $x$-pseudo-wellordering, then $\text{ot}(R) = \omega_1^x(1 + \eta) + \rho$ where $\eta = \text{ot}(\mathbb{Q})$ and $\rho < \omega_1^x$. 
\end{theorem}

\begin{proof}
See \cite{Recursive-Pseudo-Well-Orderings} or \cite{Higher-Recursion-Theory}, Lemma III.2.2.
\end{proof}

\Begin{proposition}{hyp real is pi11}
Recall if $y \in \cantorspace$, then $\text{HYP}^y = L_{\omega_1^y}(y) \cap \cantorspace$, the set of $y$-hyperarithmetic reals. 

The relation $x \in \text{HYP}^y$ is a $\Pi_1^1$ relation in the variable $x$ and $y$.
\end{proposition}

\begin{proof}
The claim is that:
$$x \in \text{HYP}^y \Leftrightarrow (\forall A)((A \in S^* \wedge y \in A) \Rightarrow (x \in A))$$
See Remark \ref{reals in omega models} about what ``$y \in A$'' should precisely mean. The latter part of the equivalence is $\Pi_1^1$. Hence the result follows from the claim.

To prove the claim: 

($\Rightarrow$) Suppose $A \in S^*$. Let $n \in \omega$ be the representative of $y$ in $A$. Since $A \models \KP$, by Lemma \ref{truncation lemma} (Truncation Lemma), $\text{WF}(A) \models \KP$. Let $\pi$ be the Mostowski collapse of $\text{WF}(A)$ onto an admissible set $B$. $y \in B$ since $y = \pi(n)$. Since $x \in \text{HYP}^y$, $x$ is in every admissible set containing $y$. $x \in B$. Then $\pi^{-1}(x)$ represents $x$ in $A$. 

($\Leftarrow$) Recall $\mathcal{HYP}(y)$ is the smallest admissible set containing $x$ and $\omega$. The domain of $\mathcal{HYP}(y)$ is $L_{\omega_1^y}(y)$. It is countable. Let $\pi : L_{\omega_1^y}(y) \rightarrow \omega$ be any bijection. The bijection gives an element $A \in S^*$ isomorphic to $\mathcal{HYP}(y)$. $\pi(y)$ represents $y$ in $A$. There exists some $n \in \omega$ such that $n$ represents $x$ in $A$, by the hypothesis. Then $x \in L_{\omega_1^y}(y)$ since $x = \pi^{-1}(n)$. $x \in \text{HYP}^y$.
\end{proof}

The following propositions uses the ideas from \cite{Higher-Recursion-Theory} III.1 and III.2.

\Begin{proposition}{tree with path no x-hyp path}
There exists a recursive tree $U$ on $2 \times \omega$ such that for all $x \in \cantorspace$, $U^x$ has a path but has no $x$-hyperarithmetic paths.
\end{proposition}

\begin{proof}
By Proposition \ref{hyp real is pi11}, there is a recursive tree $V$ on $2 \times 2 \times \omega$ such that $x \notin \text{HYP}^y$ if and only if $V^{(x,y)}$ is ill-founded. Define the relation $\Phi$ on $\cantorspace \times \bairespace$ by
$$\Phi(y,f) \Leftrightarrow (\forall n)((f_0(n) = 0 \vee f_0(n) = 1) \wedge V(f_0 \upharpoonright n, y \upharpoonright n, f_1 \upharpoonright n))$$
where $f_i(n) = f(\langle i, n \rangle)$, for $i = 0, 1$. $\Phi$ is $\Pi_1^0$. Let $U$ be a recursive tree on $2 \times \omega$ such that
$$\Phi(y,f) \Leftrightarrow (\forall n)((y \upharpoonright n, f \upharpoonright n) \in U)$$

For any $y$, if $U^y$ has a path $f$, then $\Phi(y,f)$. Therefore, $f_1 \in [V^{(f_0,y)}]$. $f_0 \notin \text{HYP}^{y}$. So $U^y$ can not have a $y$-hyperarithmetic path $f$, since otherwise $f_0 \in \text{HYP}^y$, which yields a contradiction. $U^y$ has a path: Let $x$ be any real which is not in $\text{HYP}^y$. $[V^{(x,y)}]$ is non-empty. Let $g \in [V^{(x,y)}]$. Let $f$ be such that $f_0 = x$ and $f_1 = g$. Then $\Phi(y,f)$. $f \in [U^y]$. 
\end{proof}

\Begin{definition}{kleene-brouwer ordering}
The Kleene-Brouwer ordering $<_\text{KB}$ is defined on $\finNaturalSequence$ as follows: $s <_\text{KB} t$ if and only if 

(i) $t \preceq s$ and $|t| < |s|$

or

(ii) If there exists an $n \in \omega$ such that for all $k < n$, $s(k) = t(k)$ and $s(n) < t(n)$. 
\end{definition}

\Begin{proposition}{KB ordering hyp path}
Let $T$ be a tree on $\omega$. $T$ is wellfounded if and only if $<_\text{KB}\upharpoonright T$ is wellfounded. Moreover, if there is an $x$-hyperarithmetic infinite descending sequence in $<_\text{KB} \upharpoonright T$, then there is an $x$-hyperarithmetic path through $T$.
\end{proposition}

\begin{proof}
If $f \in [T]$, then $\{f \upharpoonright n : n \in \omega\}$ is an infinite descending sequence in $<_\text{KB} \upharpoonright T$. 

Let $S = \{s_n \in \finBinarySequence : n \in \omega\}$ be an $x$-hyperarithmetic descending sequence in $<_\text{KB}\upharpoonright T$. Define $f \in \bairespace$ by
$$f(n) = i \Leftrightarrow (\exists p)(\forall q \geq p)(s_q(n) = i)$$
$f \in [T]$ and $f$ is $\Sigma_2^0(S)$. $f$ is also $x$-hyperarithmetic.
\end{proof}

Now to produce a classification of $\F$ by countable structures. The idea will be to send $x$ to an $x$-Harrison linear ordering. Using Proposition \ref{tree with path no x-hyp path} and applying the Kleene-Brouwer ordering, one can obtain a function $g$ such that $g(x)$ is an $x$-recursive $x$-pseudo-wellordering. Now suppose $\omega_1^x = \omega_1^y$. Let $\alpha$ denote this admissible ordinal. By Theorem \ref{harrison theorem}, $\text{ot}(g(x)) = \alpha(1 + \eta) + \rho_x$ and $\text{ot}(g(y)) = \alpha(1 + \eta) + \rho_y$, where $\rho_x < \alpha$ and $\rho_y < \alpha$. However, it could happen that $\rho_x \neq \rho_y$. One way to modify $g$ to get a classification of $\F$ would be to ``cut off'' the recursive tail of $g(x)$. To do this, one uses a trick, suggested Montalban, to cut off the recursive tail of the order type by taking a product of $\omega$ copies of $g(x)$. The details follow:

\Begin{proposition}{order type produce pseudo well-ordering}
Fix $x \in \cantorspace$. Let $\rho < \omega_1^x$ and $\eta = \text{ot}(\mathbb{Q})$. Then $(\omega_1^x (1 + \eta) + \rho) \omega = \omega_1^x(1 + \eta)$. 
\end{proposition}

\begin{proof}
Let $P$ be any $x$-recursive $x$-pseudo-wellorderings of order type $\omega_1^x(1 + \eta) + \rho$. Let $P \times \omega$ be the $x$-recursive structure isomorphic to $\omega$ copies of $P$ following each other. $P \times \omega$ is still an $x$-recursive $x$-pseudo-wellordering. It has no $x$-recursive tail. By Theorem \ref{harrison theorem}, $\text{ot}(P \times \omega) = (\omega_1^x(1 + \eta) + \rho)\omega = \omega_1^x(1 + \eta)$. 
\end{proof}

\Begin{proposition}{effectiveness of classification construction}
There exists an $e \in \omega$ such that for all $x \in \cantorspace$, $\{e\}^x$ is isomorphic to $(<_\text{KB} \upharpoonright U^x) \cdot \omega$, where $U$ comes from Proposition \ref{tree with path no x-hyp path}.
\end{proposition}

\begin{proof}
This is basic recursion theory using the previous results. 
\end{proof}

\Begin{theorem}{F is classifiable}
(Montalb\'{a}n) The equivalence relation $\F$ is classifiable by countable structures. In fact, there is an $e \in \omega$ such that $f(x) = \{e\}^x$ is the desired classification.
\end{theorem}

\begin{proof}
Let $\SCRL = \{\dot R\}$, where $\dot R$ is a binary relation symbol. $\F$ will be classified by countable $\SCRL$-structures. $U^x$ is an $x$-hyperarithmetic tree with paths but no $x$-hyperarithmetic path. Hence $<_\text{KB} \upharpoonright U^x$ is an $x$-recursive linear ordering with infinite descending sequences but no $x$-hyperarithmetic infinite descending sequences. So $<_\text{KB} \upharpoonright U^x$ is an $x$-recursive $x$-pseudo-wellordering. It has order type $\omega_1^x(1 + \eta) + \rho$ for some $\rho < \omega_1^x$. Therefore, $(<_\text{KB} \upharpoonright U^x) \cdot \omega$ has order type $\omega_1^x(1 + \eta)$, i.e., it is an $x$-Harrison linear ordering.  Hence $x \ \F \ y$ if and only $\omega_1^x = \omega_1^x$ if and only $\omega_1^x(1 + \eta) = \omega_1^y(1 + \eta)$ if and only $(\leq_\text{KB} \upharpoonright U^x) \cdot \omega \cong_{\SCRL} (\leq_\text{KB}\upharpoonright U^y) \cdot \omega$ if and only if $\{e\}^x \cong_{\SCRL} \{e\}^y$. This gives a classification of $\F$.
\end{proof}

\section{Finer Aspects of Classification by Countable Structures}\label{Finer Aspects of Classification by Countable Structures}

The previous section provided an explicit classification $f : \cantorspace \rightarrow S(\SCRL)$ which was $\lborel$ and for all $x \in \cantorspace$, $\text{SR}(f(x)) = \omega_1^x + 1$. This section will show that any classification of $\F$ by countable structures must have a similar property. 

The next result will calculate the complexity of each $\F$ class according to effective descriptive set theory.

\Begin{theorem}{F class not borel in any member}
For any $x \in \cantorspace$, $[x]_\F$ is not $\Pi_1^1(z)$. 
\end{theorem}

\begin{proof}
Suppose $[x]_\F$ is $\Pi^1(x)$. Let $B = \cantorspace - [x]_\F$. $B$ is then $\Sigma_1^1(x)$. Let $U$ be a tree on $2 \times \omega$ recursive in $x$ such that
$$y \in B \Leftrightarrow [U^y] \neq \emptyset$$
Observe $U \in L_{\omega_1^x}(x)$. 

Let $\SCRL$ be the language consisting of the following:

(I) A binary relation symbol $\dot\in$. 

(II) Constant symbol $\bar{a}$ for each $a \in L_{\omega_1^x}(x)$. 

(III) Two other distinguished constant symbols $\dot c$ and $\dot d$. 

\noindent $\SCRL$ can be considered a $\Delta_1$ definable subset of $L_{\omega_1^x}(x)$. 

\noindent $\SCRL$ may be regarded as a $\Delta_1$ subset of $L_{\omega_1^x}(x)$. 

Let $T$ be a theory in the countable admissible fragment $\SCRL_{L_{\omega_1^x}(x)}$ consisting of the following:

(i) $\KP$

(ii) For each $a \in L_{\omega_1^x}(x)$, $(\forall v)(v \dot \in \bar{a} \Rightarrow \bigvee_{z \in a} v = \bar{z})$. 

(iii) $\dot c \subseteq \bar{\omega}$ and $\dot d :\bar{\omega} \rightarrow \bar{\omega}$. 

(iv) For each ordinal $\sigma \in \omega_1^x$, ``$\bar{\sigma}$ is not admissible relative to $\dot c$''.

(v) $\dot d \in [\bar{U}^{\dot c}]$. 

\noindent $T$ can be considered a $\Sigma_1$ on $L_{\omega_1^x}(x)$ theory.

$T$ is consistent: Find any $y \in \cantorspace$ such that $\omega_1^y > \omega^x$. Then $y \in B$. There exists some $z \in \bairespace$ such that $z \in [U^y]$. Consider the $\SCRL$-structure $\CM$ defined as follows: $M = H_{\aleph_1}$. $\dot \in^\CM = \in \upharpoonright H_{\aleph_1}$. For each $a \in L_{\omega_1^x}(x)$, let $\bar{a}^\CM = a$. Let $\dot c^\CM = y$ and $\dot d^\CM = z$. $\CM \models T$. 

By Theorem \ref{solid model existence theorem}, $T$ has a solid model $\CN$ such that $\text{ON} \cap N = \text{ON} \cap L_{\omega_1^x}(x) = \omega_1^x$. Let $u = \dot c^\CN$ and $v = \dot d^\CM$. As in Theorem \ref{Sacks theorem}, $\omega_1^u = \omega_1^x$. 

$\CN \models v \in [U^u]$. By $\Delta_1$ absoluteness, $\text{WF}(\CN) \models v \in [U^u]$. Since $\CN$ is solid, $\text{WF}(\CN)$ is transitive as viewed in $V$. So by $\Delta_1$ absoluteness, $V \models v \in [U^u]$. $[U^u] \neq \emptyset$. $u \in B$. $\omega_1^u \neq \omega_1^x$. Contradiction.
\end{proof}

Suppose $f$ is a classification of $\F$ by countable structures in some recursive language. The Scott rank of the image of $f$ must be high:

\Begin{theorem}{classification F high scott rank}
Let $\SCRL$ be a recursive language. If $f: \cantorspace \rightarrow S(\SCRL)$ is a $\Delta_1^1(z)$ function such that $x \ \F \ y$ if and only if $f(x) \cong_{\SCRL} f(y)$, then for all $x$ such that $\omega_1^x \in \Lambda(z)$, $\text{SR}(f(x)) \geq \omega_1^x$.
\end{theorem}

\begin{proof}
Suppose there exists an $x \in \cantorspace$ with $\omega_1^x \in \Lambda(z)$ and $\text{SR}(f(x)) < \omega_1^x$. Let $\alpha = \omega_1^x$. By Proposition \ref{Sacks theorem}, there exists a $y$ with $z \leq_T y$ and $\omega_1^y = \alpha$. Since $\omega_1^y = \alpha = \omega_1^x$, $x \ \F \ y$. This implies that $f(x) \cong_{\SCRL} f(y)$. Hence $\text{SR}(f(y)) = \text{SR}(f(x)) < \omega_1^x = \alpha = \omega_1^y$. $z \leq_T y$ implies that $z \in L_{\omega_1^y}(y)$, and in particular, $z$ is in every admissible set containing $y$. Since $f$ is $\Delta_1^1(z)$, $f(y)$ is $\Delta_1^1(z,y) = \Delta_1^1(y)$ since $z \leq_T y$. $f(y)$ is hyperarithmetic in $y$. $f(y)$ is in every admissible set that has $y$ as a member. Since $\text{SR}(f(y)) < \omega_1^y$ and $f(y)$ is in every admissible set containing $y$, $\text{CSS}(f(y))$ is in every admissible set containing $y$. In particular $\text{CSS}(f(y)) \in L_{\omega_1^y}(y)$. 

By Proposition \ref{complexity of models of formula}, $\text{Mod}(\text{CSS}(f(y))$ is $\Delta_1^1(y)$. Therefore, 
$$v \in [y]_\F \Leftrightarrow f(v) \in \text{Mod}(\text{CSS}(f(y)))$$
which is $\Delta_1^1(y,z) = \Delta_1^1(y)$ since $z \leq_T y$. This contradicts Theorem \ref{F class not borel in any member}.
\end{proof}

\Begin{remark}{classification scott rank high whether highest}
Let $f$ be $\Delta_1^1(z)$ as above. For all $y \in [x]_\F$, there is an ordinal $\alpha$ such that $\text{SR}(f(y)) = \alpha$. The previous result states that if $\omega_1^x \in \Lambda(z)$, then the Scott rank of $f(x)$ is greater than or equal to $\omega_1^x$. So $\alpha \geq \omega_1^x$. Since $\omega_1^x \in \Lambda(z)$, by Theorem \ref{Sacks theorem}, there is an $x' \in \cantorspace$ such that $\omega_1^{x'} = \omega_1^x$ and $z \leq_\text{T} x'$. Then $f(x')$ is $\Delta_1^1(x', z) = \Delta_1^1(x')$. By Lemma \ref{length of scott rank}, $\text{SR}(f(x')) \leq \omega_1^{x'} + 1 = \omega_1^{x} + 1$. So one has that $\omega_1^x \leq \alpha \leq \omega_1^x + 1$. One may ask if $\alpha$ must take the largest possible value. 

Using the methods of infinitary logic as above, there is one obvious idea to try in order to force the Scott rank to be as high as possible:

Let $\mathcal{J}$ be a countable recursive language. Suppose $f : \cantorspace \rightarrow S(\mathcal{J})$ is a $\Delta_1^1(z)$ function such that $x \ \F \ y$ if and only if $f(x) \cong_\mathcal{J} f(y)$.

Since $f$ is $\Delta_1^1(z)$, it is $\Sigma_1^1(z)$. There is a tree $U$ on $2 \times 2 \times \omega$ recursive in $z$ such that $(a,b) \in f$ if and only if $[U^{(a,b)}] \neq \emptyset$. Again, one may assume $x \geq_T z$: since one can find a $x'$ with $\omega_1^{x'} = \omega_1^x$ and $x' \geq_T z$. This implies $\text{SR}(f(x')) = \text{SR}(f(x))$. 

Let $\SCRL$ be the language consisting of the following:

(I) A binary relation symbol $\dot\in$. 

(II) Constant symbols $\bar{a}$ for each $a \in L_{\omega_1^x}(x)$. 

(III) Four distinguished constant symbols $\dot c$, $\dot d$, $\dot e$, and $\dot s$.

Let $T$ be a theory in the countable admissible fragment $\SCRL_{L_{\omega_1^x}(x)}$ consisting of the following:

(i) $\mathsf{KP}$ 

(ii) For each $a \in L_{\omega_1^x}(x)$, $(\forall v)(v \dot\in \bar{a} \Rightarrow \bigvee_{u \in a} v = \bar{u})$. 

(iii) $\dot c \subseteq \bar{\omega}$, $\dot d \subseteq \bar{\omega}$, $\dot e : \bar{\omega} \rightarrow \bar{\omega}$, and $\dot s \dot\in \ {}^{<\bar{\omega}}\bar\omega$.

(iv) ``$\bar{\alpha}$ is not admissible in $\dot c$'' for each $\alpha < \omega_1^x$. 

(v) $\dot e \in [\bar{U}^{(\dot c, \dot d)}]$. 

(vi) $\rho(\dot d, \dot s) > \bar{\alpha}$ for each $\alpha < \omega_1^x$. 

\noindent $T$ can be considered a $\Sigma_1$ on $L_{\omega_1^x}(x)$ theory.

Next to show $T$ is consistent: Find $w$ such that $\omega_1^w > \omega_1^x$ and $w \in \Lambda(z)$. $(w, f(w)) \in f$, therefore, there exists some $u$ such that $u \in [U^{(w,f(w))}]$. By Theorem \ref{classification F high scott rank}, $\text{SR}(f(w)) \geq \omega_1^w$. Let $k \in {}^{<\omega}\omega$ such that $\rho(w, k) > \omega_1^x$. Define $\CM$ by $M = H_{\aleph_1}$, $\dot \in$ is the $\in$ relation of $H_{\aleph_1}$. For each $a \in L_{\omega_1^x}(x)$, $\bar{a}^\CM = a$. $\dot c^\CM = w$. $\dot d^\CM = f(w)$, and $\dot s = k$. Then $\CM \models T$. $T$ is consistent. 

By Theorem \ref{solid model existence theorem}, $T$ has a solid model $\CN$ such that $\text{ON} \cap N = \text{ON} \cap L_{\omega_1^x}(x) = \omega_1^{x}$. Let $v = \dot c^\CN$, $w = \dot d^\CN$, $u = \dot e^\CN$, and $t = \dot s^\CN$. As before, $\omega_1^v = \omega_1^x$. $\CN \models u \in [U^{(v,w)}]$. $u, v, w \in \text{WF}(N)$. By $\Delta_1$-absoluteness between transitive models, $\text{WF}(\CN) \models u \in [U^{(v,w)}]$. Since $\CN$ is solid, by $\Delta_1$-absoluteness between transitive models, $V \models u \in [U^{(v,w)}]$. Hence $w = f(v)$. 

Now, one would like to show that $\rho((w,t)) = \omega_1^x$. The problem occurs in how $\CN$ can satisfy (vi). It seems possible that there is an $\alpha < \omega_1^x$ such that for all $(v,q)$ and $\beta < \omega_1^x$, $(w,t) \sim_\alpha (v,q)$ implies $(w,t) \sim_\beta (v,q)$, but there exists some ill founded ordinal $\gamma \in N$ such that $(w,t) \not\sim_\gamma (v,q)$. That is, in $V$, $\rho((w,t)) < \omega_1^x$ but in $N$, $\rho((w,t)) > \alpha$ for all $\alpha < \omega_1^x$. 

The natural question is whether this is actually possible: Is there a structure $w$ on $\omega$, a tuple $t \in \finNaturalSequence$, and an ill-founded model $N$ of $\KP$ such that $V \models \rho((w,t)) < \text{ON} \cap N$ but for all $\alpha < \text{ON} \cap N$, $N \models \rho((w,t)) > \alpha$. 
\end{remark}

\Begin{proposition}{structure of scott rank omega10}
(Makkai) There is a hyperarithmetic (or even computable) structure $P$ such that $\text{SR}(P) = \omega_1^\emptyset$? 
\end{proposition}

\begin{proof}
See \cite{An-Example-Concerning-Scott-Heights}. Also see \cite{Computable-Trees-of-Scott-Rank}, Theorem 3.6.
\end{proof}

Before this, there had not been much difficulty proving the consistency of the desired theory by exhibiting some model with domain $H_{\aleph_1}$. A model of the next theory is not as easily produced. The classical Barwise compactness theorem will be useful in showing consistency in this case.	

\Begin{theorem}{barwise compactness}
(Barwise Compactness) Let $\CA$ be a countable admissible set and $\SCRL$ be a $\Delta_1$ in $\CA$ language. Let $\SCRL_\CA$ be the induced countable admissible fragment of $\SCRL_{\infty\omega}$. Let $T$ be a set of sentences of $\SCRL_\CA$ such that $T$ is $\Sigma_1$ in $\CA$. If every $F \subseteq T$ such that $F \in A$ has a model, then $T$ has a model.
\end{theorem}

\begin{proof}
See \cite{Admissible-Sets-and-Structures}, Theorem III.5.6. Also see \cite{Admissible-Sets}, Section 4, Corollary 8.
\end{proof}

\Begin{proposition}{rho rank bounded ordinal and unbounded illfounded ordinal}
Let $P$ be a computable structure on $\omega$ such that $\text{SR}(P) = \omega_1^\emptyset$. Then there exists an ill-founded model $N$ of $\KP$ and some $t \in \finNaturalSequence$ such that 

(i) $N \cap \text{ON} = \omega_1^\emptyset$

(ii) For all $\alpha < \omega_1^\emptyset$, $N \models \rho((P,t)) > \alpha$. 

(iii) $V \models \rho((P,t)) < \omega_1^\emptyset$. 
\end{proposition}

\begin{proof}
Let $\SCRL$ be a language consisting of the following

(I) A binary relation symbol $\dot \in$.

(II) Constant symbols $\bar{a}$ for each $a \in L_{\omega_1^\emptyset}$.

\noindent $\SCRL$ can be considered a $\Delta_1$ definable subset of $L_{\omega_1^\emptyset}$. 

Let $T$ be a theory in the countable admissible fragment $\SCRL_{L_{\omega_1^\emptyset}}$ consisting of the following

(i) $\KP$

(ii) For each $a \in L_{\omega_1^\emptyset}$, $(\forall v)(v \dot\in\bar{a} \Rightarrow \bigvee_{z \in a}v = \bar{z})$. 

(iii) For each $\alpha < \omega_1^\emptyset$, $\rho((P,\dot s)) > \bar{\alpha}$. 

\noindent $T$ can be considered a $\Sigma_1$ on $L_{\omega_1^\emptyset}$ theory. 

$T$ is consistent: Let $F \subseteq T$ such that $F \in L_{\omega_1^\emptyset}$. Then there exists $\alpha < \omega_1^\emptyset$ such that all ordinals mentioned in sentences of type (iii) are less than $\alpha$. Since $\text{SR}(P) = \omega_1^\emptyset$, there exists some $t \in \finNaturalSequence$ such that $\rho((P,t)) > \alpha$. Consider the $\SCRL$-structure defined as follows: $M = H_{\aleph_1}$. $\dot \in^\CM = \in \upharpoonright H_{\aleph_1}$. For each $a \in L_{\omega_1^\emptyset}$, $\bar{a}^\CM = a$. $\dot s^\CM = t$. Then $\CM \models F$. $F$ is consistent. By Barwise compactness (Theorem \ref{barwise compactness}), $T$ is consistent.

By Theorem \ref{solid model existence theorem}, there is a solid structure $\CN \models T$. Let $t = \dot s^\CN$. Since $\CN \models T$, for all $\alpha < \omega_1^\emptyset$, $\CN \models \rho((P,t)) > \alpha$. However, since $\text{SR}(P) = \omega_1^\emptyset$, one has $V \models \rho((P,t)) < \omega_1^\emptyset$. $\CN$ and $t$ are as desired.
\end{proof}

As mentioned before, the $\Delta_1^1$ classification $f$ of $\F$ from Theorem \ref{F is classifiable} has the property that $\text{SR}(f(x)) = \omega_1^x + 1$ for all $x$. Given the above remarks, one can ask the following

\Begin{question}{classification F high but not highest scott rank}
Does there exists a $\lborel$ function $f$ classifying $\F$ such that $\text{SR}(f(x)) = \omega_1^x$ for all $x \in \cantorspace$?
\end{question}

The authors of \cite{Computable-Trees-of-Scott-Rank} produced a very simple computable tree of Scott rank $\omega_1^\emptyset$. However, their proof in \cite{Computable-Trees-of-Scott-Rank}, Section 2 uses Barwise-Kreisel compactness and their proof in \cite{Computable-Trees-of-Scott-Rank}, Section 4 uses an overspill into the illfounded portion of the Harrison linear ordering. It is unclear if their proof method can be made uniform enough to produce in a $\borel$ manner a map taking $x$ to some $x$-relative version of their tree. 

Although it may or may not be relevant here: the distinction between structure of rank $\omega_1^x$ and $\omega_1^x + 1$ has had some role in works on the Vaught's conjecture. For example, \cite{On-the-Number-of-Countable-Models}, Theorem 4.2 shows that if $\varphi \in \SCRL_{\omega_1\omega}$ has the property that for all countable $M \models \varphi$, $\text{SR}(M) \leq \omega_1^M$, then $\varphi$ has only countably many models up to isomorphism (i.e., is not a counterexample to Vaught's conjecture).
\\*
\\*\indent Theorem \ref{classification F high scott rank} is only able to provide information about $f(x)$ when $\omega_1^x \in \Lambda(z)$ with $z$ such that $f$ is $\Delta_1^1(z)$. Some type of condition involving $\Lambda(z)$ is required: 

\Begin{lemma}{not recursivly inaccessible bijection lesser admissible}
Suppose $x \in \cantorspace$ is such that $\omega_1^x$ is not a recursively inaccessible ordinal. Then there exists a $z \in \cantorspace$ such that $z$ is $\Delta_1^1(x)$ and $\{\text{ot}(z^{[n]}) : n \in \omega\} = \Lambda(\emptyset) \cap \alpha$, where $z^{[n]} = \{y : \langle n, y \rangle \in z\}$.  
\end{lemma}

\begin{proof}
Since $\omega_1^x$ is not recursively inaccessible let $\beta$ be the largest admissible ordinal less than $\omega_1^x$. Since $\beta + 1 < \omega_1^x$, it is an $x$-recursive ordinal. There is an $e$ such that $\{e\}^{w}$ has order type $\beta + 1$. The set 
$$B = \{n \in \omega : \{e\}^x \upharpoonright n \text{ is an admissible ordinal}\}$$
 is a set in $L_{\omega_1^x}(x)$ by $\Delta_1$ separation. Let $f: \omega \rightarrow B$ be a bijection in $L_{\omega_1^x}(x)$. Now define $z$ by $z^{[n]} = \{e\}^x \upharpoonright f(n)$. 
\end{proof}

In the proof above, one needed a bijection in $L_{\alpha}(x)$ between $\omega$ and $\Lambda(\emptyset) \cap \alpha$. Note that by $\Sigma_1$ collection, there is no $\Sigma_1$ function $f : \gamma \rightarrow \alpha$ with $\gamma < \alpha$ and $f$ unbounded. If $\alpha$ is recursively inaccessible, then $\Lambda(\emptyset) \cap \alpha$ is unbounded in $\alpha$. Hence when $\alpha$ is recursively inaccessible, there can not exist such a bijection.

\Begin{proposition}{not recursively inaccessible failure below}
Suppose $\alpha < \omega_1$ is an admissible but not recursively inaccessible ordinal. Let $\SCRL = \{\dot <\}$. There exists $z$ with $\omega_1^z = \alpha$ such that 

(i) There is an $f: \cantorspace \rightarrow S(\SCRL)$ which is $\Delta_1^1(z)$.

(ii) For all $x,y \in \cantorspace$, $x \ \F \ y$ if and only if $f(x) \cong_{\SCRL} f(y)$

(iii) For all $x$ with $\omega_1^x < \omega_1^z$, $\text{SR}(f(x)) < \omega_1^x$. 
\end{proposition}

\begin{proof}
By Lemma \ref{not recursivly inaccessible bijection lesser admissible}, let $z$ be such that $\{\text{ot}(z^{[n]}) : n \in \omega\} = \Lambda(\emptyset) \cap \alpha$ and $\alpha = \omega_1^{z}$. Let $f : \cantorspace \rightarrow S(\SCRL)$ be the $\Delta_1^1$ classification given in Theorem \ref{F is classifiable}. Let $g : \omega \rightarrow S(\SCRL)$ be $\Delta_1^1$ such that for all $n \in \omega$, $g(n) \cong_\SCRL \omega + n$. 

Define the set $B \subseteq \omega \times \cantorspace$ by:
$$(m, x) \in B \Leftrightarrow \omega_1^x = \text{ot}(z^{[m]})$$
The claim is that $B$ is $\Delta_1^1(z)$:  

It is $\Sigma_1^1(z)$. 
$$(m,x) \in B \Leftrightarrow (\forall n)\Big((\{n\}^x \in \text{WO} \Rightarrow \{n\}^x \preceq z^{[m]})\Big) \wedge$$
$$(\forall k)(\exists j)(\{j\}^x \preceq z^{[m]} \wedge z^{[m]}\upharpoonright k \preceq \{j\}^x)$$

It is $\Pi_1^1(z)$. 
$$(m, x) \notin B \Leftrightarrow (\exists j)(\forall n)(\{n\}^x \in \text{WO} \Rightarrow (\{n\}^x \preceq z^{[m]} \upharpoonright j) \vee (\exists n)(\{n\}^x \in \text{LO} \wedge \{n\}^x \preceq z^{[m]} \wedge z^{[m]} \preceq \{n\}^x)$$

Now define the following function $h : \cantorspace \rightarrow S(\SCRL)$. 

$$(x,y) \in h \Leftrightarrow (\exists n)\Big((n,x) \in B \wedge y = g(n)\Big) \vee (\forall n)\Big((n,x) \notin B \wedge y = f(x)\Big)$$

$f$ is $\Delta_1^1(z)$. For all $x,y$, $x \ \F \ y$ if and only if $f(x) \cong_\SCRL f(y)$. If $\omega_1^x < \omega_1^z$, then $\text{SR}(f(x)) = \text{SR}(\omega + n) < \omega_1^x$, where $n$ is such that $\text{ot}(z^{[n]}) = \omega_1^x$. 
\end{proof}

Proposition \ref{not recursively inaccessible failure below} asserts that for each $\alpha < \omega_1$ which is admissible but not recursively inaccessible, there exists some $z$ with $\omega_1^z = \alpha$ and some $\Delta_1^1(z)$ classification of $\F$ such that the Scott rank condition of Theorem \ref{classification F high scott rank} fails on all the $\F$-classes associated with admissible ordinals less than $\alpha$. Can this also be achieved when $\alpha$ is recursively inaccessible?

The most interesting question of this kind is: Is there some classification $f$ of $\F$ which is $\Delta_1^1(z)$ and the Scott rank condition fails for some class associated with an admissible ordinal $\alpha > \omega_1^z$? 

\section{Almost Borel Reductions}\label{Almost Borel Reductions}
\Begin{definition}{pinned cardinal}
(\cite{Forcing-Borel-Reducibily-Invariants} Definition 3.1.1) Let $E$ be a $\bSigma_1^1$ equivalence relation on a Polish space $X$. Let $\bbP$ be a forcing and $\tau$ be a $\bbP$-name for an element of $X$, i.e. $1_\bbP \forces_\bbP \tau \in X$. Let $\tau_\text{left}$ and $\tau_\text{right}$ be $\bbP^2$-names for the evaulation of $\tau$ according to the left and right $\bbP$-generic coming from the $\bbP^2$-generic. $\tau$ is an $E$-pinned $\bbP$-name if and only if $1_{\bbP^2} \forces \tau_\text{left} \ E \ \tau_\text{right}$.

(\cite{Forcing-Borel-Reducibily-Invariants} Definition 4.1.1 and 4.1.2) Let $\bbP$ and $\bbQ$ be two forcings. $\sigma$ be an $E$-pinned $\bbP$ name and $\tau$ be an $E$-pinned $\bbQ$-name. Define the relation $\sigma \ \bar{E} \ \tau$ if and only if $\bbP\times\bbQ \forces \sigma \ E \ \tau$ (where $\sigma$ and $\tau$ are considered $\bbP \times \bbQ$-names in the natural way). The pinned cardinal of $E$, denoted $\kappa(E)$, is the smallest cardinal $\kappa$ such that every $E$-pinned $\bbP$-name is $\bar{E}$-related to an $E$-pinned $\bbQ$-name with $|\bbQ| < \kappa$, if this cardinal exists. Otherwise, $\kappa(E) = \infty$.
\end{definition}

\Begin{definition}{Ew1}
$E_{\omega_1}$ is the $\Sigma_1^1$ equivalence relation on $\cantorspace$ defined by $x \ E_{\omega_1} \ y$ if and only if $(x \notin \text{WO} \wedge y \notin \text{WO}) \vee (\text{ot}(x) = \text{ot}(y))$. 
\end{definition}

\Begin{proposition}{Ew1 pinned cardina}
$\kappa(E_{\omega_1}) = \infty$
\end{proposition}

\begin{proof}
See \cite{Forcing-Borel-Reducibily-Invariants}, Example 4.1.8.
\end{proof}

\Begin{definition}{aleph0 almost borel reduciblity}
Let $E$ and $F$ be two equivalence relations on Polish spaces $X$ and $Y$, respectively. $E \leq_{\text{a}\borel} Y$ if and only if there is a $\bDelta_1^1$ function $f : X \rightarrow Y$ and a countable set $A \subseteq X$ such that if $c, d \notin A$, then $c \ E \ d$ if and only if $f(c) \ F \ f(d)$. In this situation, one says $E$ is almost $\bDelta_1^1$ reducible to $F$. (It is called a weak Borel reduction in \cite{Forcing-Borel-Reducibily-Invariants} Definition 2.1.2.)
\end{definition}

\Begin{proposition}{pinned cardinal aleph0 almost borel reducibility}
Let $E$ and $F$ be $\bSigma_1^1$ equivalence relations on Polish spaces $X$ and $Y$, respectively. If $E \leq_{\text{a}\borel} F$, then $\kappa(E) \leq \kappa(F)$. 
\end{proposition}

\begin{proof}
See \cite{Forcing-Borel-Reducibily-Invariants}, Theorem 4.1.3.
\end{proof}

\Begin{proposition}{F pinned cardinal}
$\kappa(\F) = \infty$.
\end{proposition}

\begin{proof}
For any cardinal $\kappa$, consider the forcing $\text{Coll}(\omega, \kappa)$. Let $\tau$ be a $\text{Coll}(\omega, \kappa)$ name for a real such that $1_{\text{Coll}(\omega, \kappa)} \forces_{\text{Coll}(\omega, \kappa)} \omega_1^\tau = \check\kappa$. $\tau$ is a $\F$-pinned $\text{Coll}(\omega, \kappa)$-name, since 
$$1_{\text{Coll}(\omega,\kappa) \times \text{Coll}(\omega,\kappa)} \forces_{\text{Coll}(\omega,\kappa) \times \text{Coll}(\omega, \kappa)} \omega_1^{\tau_\text{left}} = \check \kappa = \omega_1^{\tau_\text{right}}$$
Now suppose $\bbQ$ is a forcing and $\sigma$ is an $\F$-pinned $\bbQ$-name with $\tau \ \bar{\F} \ \sigma$. This implies that $1_\bbQ \forces_{\bbQ} \omega_1^\sigma = \check \kappa$. $1_\bbQ \forces_\bbQ |\check\kappa| = \aleph_0$. Since any forcing $\bbQ$ is $|\bbQ|^+$-cc. $\bbQ$ preserves cardinals greater than or equal to $|\Q|^+$. Since $\bbQ$ makes $\kappa$ countable, $|\bbQ| \geq \kappa$. $\kappa(\F) \geq \kappa$. Since $\kappa$ was arbitrary, $\kappa(\F) = \infty$. 
\end{proof}

\Begin{theorem}{Ew1 aleph0 almost reduction pinned cardinal}
(Zapletal) Suppose there exists a measurable cardinal. Let $E$ be a $\bSigma_1^1$ equivalence relation. $\kappa(E) = \infty$ if and only if $E_{\omega_1} \leq_{\text{a}\borel} E$. 
\end{theorem}

\begin{proof}
See \cite{Forcing-Borel-Reducibily-Invariants}, Theorem 4.2.1.
\end{proof}

\Begin{proposition}{measurable F aleph1 almost reduction}
($\mathsf{ZFC + Measurable\ Cardinal}$) $E_{\omega_1} \leq_{\text{a}\borel} \F$. 
\end{proposition}

\begin{proof}
This follows from Theorem \ref{Ew1 aleph0 almost reduction pinned cardinal} and Proposition \ref{F pinned cardinal}.
\end{proof}

Since Theorem \ref{Ew1 aleph0 almost reduction pinned cardinal} assumes a measurable cardinal, a natural task would be to investigate the consistency strength of the statement ``For all $\bSigma_1^1$ equivalence relation, $\kappa(E)= \infty$ if and only if $E_{\omega_1} \leq_{\text{a}\borel} E$''.

Therefore, an interesting question is whether $L$ satisfies the above statement. The rest of this section will consider this question.

\Begin{theorem}{complexity of Ew1 classes}
Suppose $x \in \text{WO}$ and $y \in \cantorspace$ such that $\omega_1^y < \text{ot}(x)$, then $[x]_{E_{\omega_1}}$ is not $\Sigma_1^1(y)$.
\end{theorem}

\begin{proof}
Suppose $[x]_{E_{\omega_1}}$ was $\Sigma_1^1(y)$. Let $U$ be a tree on $2 \times \omega$ which is recursive in $y$ and 
$$(\forall u)(u \in [x]_{E_{\omega_1}} \Leftrightarrow (\exists f)(f \in [U^u]))$$

\noindent Let $\SCRL$ be a language consisting of the following:

(i) A binary relation symbol $\dot\in$.

(ii) For each $a \in L_{\omega_1^y}(y)$, a constant symbol $\bar{a}$.

(iii) Two distinct constant symbols $\dot c$ and $\dot d$.

\noindent $\SCRL$ may be considered a $\Delta_1$ definable language over $L_{\omega_1^y}(y)$.

Let $T$ be a theory in the countable admissible fragment $\SCRL_{L_{\omega_1^y}(y)}$ consisting of the following sentences:

(I) $\mathsf{KP}$ 

(II) For each $a \in L_{\omega_1^y}(y)$, $(\forall v)(v \dot\in \bar{a} \Rightarrow \bigvee_{z \in a} v = \bar{z})$.

(III) $\dot c \subseteq \bar{\omega}$, $\dot d : \bar{\omega} \rightarrow \bar\omega$

(IV) $\dot d \in [U^{\dot c}]$.

(V) For all $\alpha < \omega_1^y$, $\bar{\alpha}$ is not admissible in $\dot c$.

\noindent $T$ may be considered a $\Sigma_1$ theory in $L_{\omega_1^y}(y)$.

Next, the claim is that $T$ is consistent. Since $x \in [x]_{E_{\omega_1}}$, there exists $g$ such that $g \in [U^{x}]$. Consider the $\SCRL$-structure $\CN$ defined as follows: Let the universe $N$ be $H_{\aleph_1}$. Let $\dot\in^\CN = \in \upharpoonright H_{\aleph_1}$. Let $\dot c^\CN = x$ and $\dot d^\CN = g$. $\CN \models T$. For (V), observe that if $\CA$ is an admissible set with $x \in A$, then $\text{ot}(x) \in A$. Hence $\text{ON} \cap A > \text{ot}(x) > \omega_1^y$.

By Theorem \ref{solid model existence theorem}, let $\CM$ be a solid model of $T$. Let $z = c^\CM$. $z \in [x]_{E_{\omega_1}}$ since $\dot d^\CM \in [U^z]$. As in the proof of Sacks theorem, $\omega_1^z = \omega_1^y$. $z \in L_{\omega_1^z}(z)$. So $\text{ot}(z) \in L_{\omega_1^z}(z)$. This is impossible since $\omega_1^z = \omega_1^y < \text{ot}(x) = \text{ot}(z)$.
\end{proof}

\Begin{theorem}{spectrum of aleph0 borel reduction Ew1 to F}
If $f: \cantorspace \rightarrow \cantorspace$ is $\Delta_1^1(y)$ and witnesses $E_{\omega_1} \leq_{\text{a}\borel} \F$, then there exists a $\beta < \omega_1$ such that for all $\alpha \in \Lambda(y)$ with $\alpha > \beta$, the next admissible ordinal after $\alpha$ is not in $\Lambda(y)$.
\end{theorem}

\begin{proof}
Let $f : \cantorspace \rightarrow \cantorspace$ witness $E_{\omega_1} \leq_{\text{a}\borel} \F$. There exists some countable set $A \subseteq \cantorspace$ such that $x \ E_{\omega_1} \ y$ if and only if $f(x) \ \F \ f(y)$ whenever $x,y \notin A$. Let $\beta = \sup \{\text{ot}(x) : x \in A\}$. The claim is that this $\beta$ works. So suppose not. There exists $\alpha', \alpha \in \Lambda(y)$ such that $\alpha > \beta$, $\alpha' > \beta$, and $\alpha$ is the next admissible ordinal after $\alpha'$.

Since $f$ is $\Delta_1^1(y)$, let $U$ be a tree on $2 \times 2\times \omega$ such that for all $a,b \in \cantorspace$, $(a,b) \in f \Leftrightarrow [U^{(a,b)}]$ is ill-founded.

Claim: There exists $a, b \in \cantorspace$ such that $\alpha' < \text{ot}(a) < \text{ot}(b) < \alpha$, $\omega_1^{f(a)} \geq \alpha$, and $\omega_1^{f(b)} \geq \alpha$.

To prove this claim: 
If there exists a $c \in \cantorspace$ such that $\alpha' < \text{ot}(c) < \alpha$ and $\omega_1^{f(c)} = \alpha'$, then fix such a $c$. If not, pick any $c \in \cantorspace$ such that $\alpha' < \text{ot}(c) < \alpha$. In this latter case, $c$ will just be ignored.

Then any $d \in \cantorspace$ with $\text{ot}(d) > \beta$, $d \notin [c]_{E_{\omega_1}}$, $\omega_1^{f(d)} \neq \alpha'$ since $f$ is a reduction. Pick any $d \in \cantorspace$ with $d \notin [c]_{E_{\omega_1}}$ and $\alpha' < \text{ot}(d) < \alpha$.

Suppose that $\omega_1^{f(d)} < \alpha'$. By Proposition \ref{Sacks theorem}, let $z$ be any real such that $\omega_1^z = \alpha'$ and $y \leq_T z$. $[f(d)]_\F$ is $\Delta_1^1(z)$ by Proposition \ref{F classes borel in real of higher church kleene}.
$$k \in [d]_{E_{\omega_1}} \Leftrightarrow f(k) \in [f(d)]_\F$$
Hence $[d]_{E_{\omega_1}}$ is $\Sigma_1^1(y, z) = \Sigma_1^1(z)$. However, $\omega_1^z = \alpha' < \text{ot}(d)$. This contradicts Theorem \ref{complexity of Ew1 classes}. 

This shows that $\omega_1^{f(d)} \geq \alpha'$. Since $d \notin [c]_{E_{\omega_1}}$, $\omega_1^{f(d)} > \alpha'$. However, the next admissible ordinal greater than $\alpha'$ is $\alpha$. Therefore, $\omega_1^{f(d)} \geq \alpha$.

Now let $a, b$ be any two reals such that $a,b \notin [c]_{E_{\omega_1}}$ and $\alpha' < \text{ot}(a) < \text{ot}(b) < \alpha$. Since $d$ in the above was arbitrary with these two properties, these two reals satisfy Claim.

Now fix $a,b \in \cantorspace$ satisfying the claim. $\text{SR}(a) = \text{ot}(a)$ and $\text{SR}(b) = \text{ot}(b)$. Thus their canonical Scott sentence are both elements of $L_\alpha$ since $\text{ot}(a), \text{ot}(b) \in L_\alpha$, $\text{SR}(\text{ot}(a)) < \alpha$, $\text{SR}(\text{ot}(b)) < \alpha$, and Proposition \ref{phi formula}.

Let $\SCRL$ be a language consisting of:

(i) A binary relation symbol $\dot\in$.

(ii) For each $e \in L_{\alpha}(y)$, a constant symbol $\bar{e}$.

(iii) Six distinct symbols $\dot a, \dot b, \dot c, \dot d, \dot u, \dot v$.

\noindent $\SCRL$ may be considered as a $\Delta_1$ definable language in $L_\alpha(y)$.

Let $T$ be a theory in the countable admissible fragment $\SCRL_{L_{\alpha}(y)}$ consisting of the following sentences:

(I) $\mathsf{KP}$ in the symbol $\dot\in$.

(II) For each $e \in L_\alpha(y)$, $(\forall v)(v \dot\in \bar{e} \Rightarrow \bigvee_{z \in e}v = \bar{z})$.

(III) $\dot a, \dot b, \dot c, \dot d \subseteq \bar{\omega}$. $\dot u, \dot v$ are functions from $\bar{\omega} \rightarrow \bar{\omega}$.

(IV) $\dot u \in [U^{(\dot a, \dot c)}]$ and $\dot v \in [U^{(\dot b, \dot d)}]$. 

(V) $\dot a \models \text{CSS}(a)$ and $\dot b \models \text{CSS}(b)$.

(VI) For all $\beta < \alpha$, $\bar\beta$ is not admissible in $\dot c$ and $\bar{\beta}$ is not admissible in $\dot d$.

\noindent $T$ may be considered a $\Sigma_1$ theory in $L_\alpha(y)$.

Since $(a, f(a)) \in f$ and $(b, f(b)) \in f$, let $u, v \in \bairespace$, be such that $u \in [U^{(a,f(a))}]$ and $v \in [U^{(b,f(b))}]$. 

To show to that $T$ is consistent: consider the following model of $\CN$: The universe $N$ is $H_{\aleph_1}$. $\dot\in^\CN = \in \upharpoonright H_{\aleph_1}$. For each $e \in L_\alpha$, $\bar{e}^\CN = e$. $\dot a^\CN = a$, $\dot b^\CN = b$. Let $\dot c^\CN = f(a)$ and $\dot d^\CN = f(b)$. Let $\dot u^\CN = u$ and $\dot v^\CN = v$. Then $\CN \models T$. 

By Theorem \ref{solid model existence theorem}, there exists a solid model $\CM \models T$. Let $a' = \dot a^\CM$, $b' = \dot b^\CM$. $f(a') = \dot c^\CM$ and $f(b') = \dot d^\CM$ since $\dot u^\CM \in [U^{(a', \dot c^\CM)}]$ and $\dot v^\CM \in [U^{(b', \dot d^{\CM})}]$. As in the proof of Sacks theorem, $\omega_1^{f(a')} = \omega_1^{f(b')} = \alpha$. By absoluteness of satisfaction from $\CM$ to $\text{WF}(\CM)$ to $V$, $a' \models \text{CSS}(a)$ and $b' \models \text{CSS}(b)$. Hence in $V$, $\text{ot}(a') = \text{ot}(a)$ and $\text{ot}(b') = \text{ot}(b)$. In particular, $\text{ot}(a') \neq \text{ot}(b')$. Hence $\neg (a' \ E_{\omega_1} \ b')$. However, $\omega_1^{f(a')} = \omega_1^{f(b')} = \alpha$ implies $f(x) \ \F \ f(y)$. This contradicts $f$ being a reduction.

This proves the theorem for those $\alpha \in \Lambda(y) \cap \omega_1$. Note the statement that $f$ and countable $A \subseteq \cantorspace$ witnesses $E_{\omega_1} \leq_{\text{a}\borel} \F$ can be written as
$$(\forall x)(\forall y)((x \notin A \wedge y \notin A) \Rightarrow (x \ E_{\omega_1} \ y \Leftrightarrow f(x) \ \F \ f(y)))$$
This is $\Pi_2^1(y,A)$ and so holds in all generic extensions by Schoenfield's absoluteness. To show the theorem holds for all $\alpha \in \Lambda(y)$ and $\alpha \geq \omega_1$, let $G \subseteq \text{Coll}(\omega, \alpha)$ be $\text{Coll}(\omega, \alpha)$-generic over $V$. In $V[G]$, let $\beta = \sup\{\text{ot}(x) : x \in A\}$ be the same ordinal as before. Since $\beta < \omega_1^V \leq \alpha$, the result above, applied in $V[G]$ for $\Lambda(y) \cap \omega_1^{V[G]}$, will show the theorem holds for $\alpha$. This concludes the proof.
\end{proof}

\Begin{theorem}{no almost reduction Ew1 F in L}
$L \models \neg(E_{\omega_1} \leq_{\text{a}\borel} \F)$. This also holds in set generic extensions of $L$. 
\end{theorem}

\begin{proof}
In $L$, for all $x \in \cantorspace$, there exists some $\alpha < \omega_1$ such that $x \in L_\alpha$. Then $\Lambda(x) - \alpha = \Lambda(\emptyset) - \alpha$. Hence there are no reals with admissible spectrum as described in Theorem \ref{spectrum of aleph0 borel reduction Ew1 to F}.
\end{proof}

\section{Counterexamples to Vaught's Conjecture and $\F$}\label{counterexamples}

\Begin{definition}{isomorphism equiv of models}
Let $\SCRL$ be a recursive language. Let $\varphi \in \SCRL_{\omega_1\omega}$. Define $E_\SCRL^\varphi$ to be the $\analytic$ equivalence relation on $S(\SCRL)$ defined by
$$x \ E_\SCRL^\varphi \ y \Leftrightarrow (x \not\models \varphi \wedge y \not\models \varphi) \vee (x \cong_\SCRL y)$$
\end{definition}

See Proposition \ref{complexity of models of formula} for the $\borel$ definability of $x \not\models \varphi$. 

\Begin{definition}{def counterexample vaught conjecture}
A \textit{counterexample to Vaught's conjecture} is a $\varphi \in \SCRL_{\omega_1\omega}$ (for some recursive language $\SCRL$) such that $E_\SCRL^\varphi$ is a thin equivalence relation with uncountably many classes.
\end{definition}

From a list of questions from the Vaught's Conjecture Workshop 2015 at the University of California at Berkeley, Sy-David Friedman asked the following question:

\Begin{question}{sy friedman vaught question}
(Sy-David Friedman) Is there some recursive language $\SCRL$ such that $\F$ is $\borel$ bireducible to the $\SCRL$-isomorphism relation restricted to some $\borel$ invariant set? 
\end{question}

Every invariant $\borel$ set for the $\SCRL$-isomorphism relation is of the form $\text{Mod}(\varphi)$ for some $\varphi \in \SCRL_{\omega_1\omega}$ (see \cite{Invariant-Descriptive-Set-Theory}, Theorem 11.3.6). Therefore, the above question is equivalent to whether there exists some $\varphi$ such that $\F \equiv_\borel E_\SCRL^\varphi$. 

$E_\SCRL^\varphi \leq_\borel \F$ implies that $E_\SCRL^\varphi$ is thin. $\F \leq_\borel E_\SCRL$ implies that $E_{\SCRL}^\varphi$ has uncountably many classes. Hence any such $\varphi$ is a counterexample to Vaught's conjecture. 

Using the ideas from the previous section, it will be shown that in $L$, no counterexample $\varphi$ of Vaught's conjecture has the property that $E_\SCRL^\varphi \leq_\borel \F$. Hence, Friedman's question has a negative answer in $L$. 

\Begin{theorem}{isomorphism class not analytic in some real of small omega1}

Let $\SCRL$ be a recursive language. Let $M \in S(\SCRL)$ and $y \in \cantorspace$ be such that $\omega_1^y < \text{R}(M)$. Then $[M]_{\cong_\SCRL}$ is not $\Sigma_1^1(y)$. ($\cong_\SCRL$ denote the equivalence relation of $\SCRL$-isomorphism. Recall $R$ is defined in Definition \ref{scott relation and sentence}.)
\end{theorem}

\begin{proof}
Suppose $[M]_{\cong_\SCRL}$ is $\Sigma_1^1(y)$. Let $U$ be a tree on $2 \times \omega$ which is recursive in $y$ and 
$$(\forall N)(N \in [M]_{\cong_\SCRL} \Leftrightarrow (\exists f)(f \in [U^N]))$$
Let $\mathscr{U}$ be the language consisting of the following:

(i) A binary relation symbol $\dot\in$. 

(ii) For each $a \in L_{\omega_1^y}(y)$, a constant symbol $\bar{a}$. 

(iii) Two distinct constant symbols $\dot c$ and $\dot d$. 

\noindent $\mathscr{U}$ may be considered a $\Delta_1$ definable class in $L_{\omega_1^y}(y)$. 

Let $T$ be a theory in the countable admissible fragment $\mathscr{U}_{L_{\omega_1^y(y)}}$ consisting of the following sentences:

(I) $\KP$

(II) For each $a \in L_{\omega_1^y}(y)$, $(\forall v)(v \dot\in \bar{a} \Rightarrow \bigvee_{z \in a} v = \bar{z})$. 

(III) $\dot c \subseteq \bar\omega$, $\dot d : \bar{\omega} \rightarrow \bar{\omega}$. 

(IV) $\dot d \in [U^{\dot c}]$. 

(V) For all $\alpha < \omega_1^y$, $\bar\alpha$ is not admissible in $\dot c$. 

\noindent $T$ may be considered a $\Sigma_1$ definable theory in $L_{\omega_1^y}(y)$. 

$T$ is consistent: Since $M \in [M]_{\cong_\SCRL}$, there is some $g$ such that $g \in [U^M]$. Define a $\SCRU$-structure $\CN$ as follows: Let the universe $N$ be $H_{\aleph_1}$. Let $\dot\in^\CN = \dot\in\upharpoonright H_{\aleph_1}$. Let $\dot c^\CN = M$ and $\dot d^\CN = g$. $\CN$ is a model of $T$. To see (V), note that $\omega_1^M > \omega_1^y$. This is because if otherwise, $M$ would be an element of some admissible set $\CA$ such that $o(\CA) \leq \omega_1^y$. By Proposition \ref{length of scott rank}, $\text{R}(M) \leq \text{ON} \cap A \leq \omega_1^y$, which is a contradiction.

By Theorem \ref{solid model existence theorem}, let $\CM$ be a solid model of $T$. Let $P = \dot c^\CM$. $P \in [M]_{\cong_\SCRL}$ since $\dot d \in [U^{P}]$. Like in the proof of Sacks theorem, $\omega_1^P = \omega_1^y$. Therefore, $P \in L_{\omega_1^y}(P)$. By Proposition \ref{length of scott rank}, $\text{R}(P) \leq \omega_1^y$. However $P \in [M]_{\cong_\SCRL}$ implies that $P \cong_{\SCRL} M$. $\text{R}(P) \leq \omega_1^y < R(M)$. Contradiction. 
\end{proof}

\Begin{fact}{counterexample vaught model of limit rank}
Let $\SCRL$ be a recursive language. If $\varphi$ is a counterexample to the Vaught conjecture, then for all limit ordinals $\beta > \text{qr}(\varphi)$, $\varphi$ has a model of Scott rank $\beta$. 
\end{fact}

\begin{proof}
See \cite{Scott-Processes}, Theorem 10.8. A similar result is also shown in \cite{Scott-Rank-of-Counterexamples-to-Vaught's-Conjecture}, page 19. 
\end{proof}

\Begin{fact}{real coding ordinal in certain admissible sets}
Let $\CA$ be a countable admissible set. Let $\alpha < o(\CA)$. Then there exists a countable admissible set $\CB$ extending $\CA$ with $o(\CB) = o(\CA)$ and such that there exists a $c \in \cantorspace \cap B$ with $\text{ot}(c) = \alpha$. 
\end{fact}

\begin{proof}
This can be proved using the techniques of infinitary logic in the countable admissible fragment $\CA$ using a Scott sentence for $\alpha$ as a linear ordering. Since this is similar to several previous arguments, the details are omitted. 
\end{proof}

\Begin{fact}{harrington silver theorem}
Let $E$ be a $\Pi_1^1(z)$ equivalence relation on a Polish space $\cantorspace$ with countably many classes. Then for all $x \in \cantorspace$, there is a $\lborel(z)$ set $U$ such that $x \in U \subseteq [x]_E$. 
\end{fact}

\begin{proof}
See \cite{A-Powerless-Proof-of-a-Theorem-of-Silver}. 

In the effective proof of Silver's dichotomy for $\bPi_1^1$ equivalence using the Gandy-Harrington topology, the two outcomes depend on whether the set
$$V = \{x \in \cantorspace : \text{There exists $\lborel(z)$ set $U$ with } x \in U \subseteq [x]_E\}$$
is equal to $\cantorspace$. If $V = \cantorspace$, then $E$ has only countable many classes. This gives the desired result above. See \cite{Invariant-Descriptive-Set-Theory}, Theorem 5.3.5 for a presentation of the effective proof of Silver's theorem.
\end{proof}

\Begin{fact}{model of high omega1}
Let $\SCRL$ be a recursive language. Let $\varphi \in \SCRL_{\omega_1\omega}$ be a counterexample to Vaught's conjecture. Let $z \in \cantorspace$ be such that $\varphi \in L_{\omega_1^z}(z)$. Let $\beta$ be a $z$-admissible ordinal. Let $M, N \in S(\SCRL)$ be such that $M \models \varphi$, $N \models \varphi$, $R(M) < \beta$ and $R(N) < \beta$. Then there exists a countable admissible set $\CA$ extending $L_{\beta}(z)$ such that $o(\CA) = \beta$ and $\text{CSS}(M), \text{CSS}(N) \in A$. 
\end{fact}

\begin{proof}
Let $\alpha < \beta$ be such that $\text{R}(M) < \alpha$ and $\text{R}(N) < \alpha$. $L_{\beta}(z)$ is an admissible set. By Fact \ref{real coding ordinal in certain admissible sets}, there exists some countable admissible set $\CB$ extending $L_\beta(z)$ containing some real $c$ which codes $\alpha$ and $o(\CB) = \beta$. 

Let $\equiv_\alpha^\SCRL$ be the relation of $\SCRL$-elementary equivalence with respect to just the formulas of quantifier rank less than $\alpha$. $\equiv_\alpha^\SCRL \upharpoonright \text{Mod}(\varphi)$ is a $\lborel(c,z)$ equivalence relation. (Proposition \ref{complexity of models of formula} is used here.) Since $\varphi$ is a counterexample to Vaught's conjecture, $\equiv_\alpha^\SCRL \upharpoonright \text{Mod}(\varphi)$ has only countably many classes. By Fact \ref{harrington silver theorem}, there exists $\lborel(c,z)$ sets $U_M$ and $U_N$ such that $M \in U_M \subseteq [M]_{\equiv_\alpha^\SCRL \upharpoonright \text{Mod}(\varphi)}$ and $N \in U_N \subseteq [N]_{\equiv_\alpha^\SCRL\upharpoonright\text{Mod}(\varphi)}$. Let $T_M$ and $T_N$ be the $c\oplus z$ recursive trees such that for all $P$
$$P \in U_M \Leftrightarrow T_M^P \text{ is illfounded}$$
$$P \in U_N \Leftrightarrow T_N^P \text{ is illfounded}$$

Let $\SCRU$ be the language consisting of the following:

(i) A binary relation symbol $\dot\in$.

(ii) For each $a \in B$, a constant symbol $\bar{a}$. 

(iii) Four new constant symbols, $\dot R$, $\dot S$, $\dot e$, and $\dot f$. 

\noindent $\SCRU$ may be considered a $\Delta_1$ definable language in $\CB$. 

Let $T$ be a theory in the countable admissible fragmant $\SCRU_{\CB}$ consisting of the following sentences:

(I) $\KP$

(II) For each $a \in B$, $(\forall v)(v \dot \in \bar{a} \Rightarrow \bigvee_{z \in a} v = \bar{z})$.

(III) $\dot R \subseteq \bar{\omega}$, $\dot S \subseteq \bar{\omega}$, $\dot e \in [T_M^{\dot R}]$, and $\dot f \in [T_N^{\dot R}]$. 

\noindent $T$ may be considered as a $\Sigma_1$ definable set in $\CA$. 

$T$ is consistent: Since $M \in U_M$ and $N \in U_N$, find some $v \in [T_M^M]$ and $w \in [T_N^N]$. Consider the $\SCRU$ structure $\mathcal{I}$ defined as follows: Let the universe $I$ be $H_{\aleph_1}$. Let $\dot\in^\mathcal{I} = \in \upharpoonright H_{\aleph_1}$. Let $\dot R^\mathcal{I} = M$, $\dot S^\mathcal{I} = N$, $\dot e^\mathcal{I} = v$, and $\dot f^\mathcal{I} = w$. Then $\mathcal{I} \models T$. 

By Theorem \ref{solid model existence theorem}, let $\mathcal{J}$ be a solid model of $T$ with $o(\mathcal{J}) = o(\CB)$. Let $R = \dot R^\mathcal{J}$, $S = \dot S^\mathcal{J}$, $e = \dot e^\mathcal{J}$, and $f = \dot f^\mathcal{J}$. By $\Delta_1$ absoluteness (first between $\mathcal{J}$ and $\text{WF}(\mathcal{J})$ and then between $\text{WF}(\mathcal{J})$ and $V$), $e \in [T_M^R]$ and $f \in [T_M^S]$. Hence $R \in U_M$ and $S \in U_N$. 

Let $\CA = \text{WF}(\mathcal{J})$. By Lemma \ref{truncation lemma}, $\CA$ is an admissible set. It has been shown that $\CA$ has two elements $R$ and $S$ such that $R \in [M]_{\equiv_\alpha^\SCRL \upharpoonright \text{mod}(\varphi)}$ and $S \in [N]_{\equiv_\alpha^\SCRL \upharpoonright \text{Mod}(\varphi)}$. 

Since $\text{CSS}(M)$ and $\text{CSS}(N)$ has quantifier rank less than $\alpha$, $M \equiv_\alpha^\SCRL R$, and $N \equiv_\alpha^\SCRL$, the following must hold: $R \models \text{CSS}(M)$ and $S \models \text{CSS}(N)$. Hence $\text{CSS}(R) = \text{CSS}(M)$ and $\text{CSS}(S) = \text{CSS}(N)$. 

Since $R(M), R(N) < \alpha$, Proposition \ref{phi formula} implies that $\text{CSS}(R) \in A$ and $\text{CSS}(S) \in A$. Therefore, $\text{CSS}(M) \in A$ and $\text{CSS}(N) \in A$. This completes the proof.
\end{proof}

\Begin{theorem}{counterexample vaught not reducible F}
Let $\SCRL$ be a recursive language. Let $\varphi \in \SCRL_{\omega_1\omega}$ be a counterexample to Vaught's conjecture. Suppose $f$ is a $\borel$ function witnessing $E_\SCRL^\varphi \leq_\borel \F$, then there exists some ordinal $\gamma$ and real $z$ such that for all $\alpha \in \Lambda(z)$ with $\alpha > \gamma$, the next admissible ordinal greater than $\alpha$ is not in $\Lambda(z)$. 
\end{theorem}

\begin{proof}
First, the theorem will be shown for $\Lambda(z) \cap \omega_1$. At the end, this result will be used to obtain the theorem for the full $\Lambda(z)$. 

Let $f : S(\SCRL) \rightarrow \cantorspace$ be $\lborel(r)$ witnessing $E_\SCRL^\varphi \leq_\borel \F$ where $r$ is some real. Find any $s \in \cantorspace$ such that $\varphi \in L_{\omega_1^s}(s)$. Let $z = r \oplus s$. Note that $f$ is $\lborel(z)$. Let $\gamma = \omega_1^z$. Certainly $\gamma > \text{qr}(\varphi)$. 

Now suppose there exists some $\alpha, \beta \in \Lambda(z) \cap \omega_1$ with $\gamma < \alpha$ and $\beta$ is the next admissible ordinal greater than $\alpha$. 

Between two consecutive admissible ordinals, there are infinitely many limit ordinals. Since $\varphi$ is a counterexample to the Vaught's conjecture, Fact \ref{counterexample vaught model of limit rank} implies that there are infinitely many models of $\varphi$ with Scott ranks between $\alpha$ and $\beta$. Let $P$, $M$, and $N$ be three models of $\varphi$ with distinct Scott rank between $\alpha$ and $\beta$. Since $f$ is a reduction of $E_\SCRL^\varphi$ to $\F$, at most one $X \in \{P, M, N\}$ has the property that $\omega_1^{f(X)} = \alpha$. If such an $X$ among these three exists, then without loss genererality, assume it was $P$. (If no $X$ among these three has this property, then one can just ignore $P$ for the rest of the proof.) 

Now to show that $\omega_1^{f(M)} \geq \beta$ and $\omega_1^{f(N)} \geq \beta$: Suppose $\omega_1^{f(M)} < \beta$. Since $P$ and $M$ are not $\SCRL$-isomorphic and $f$ is a reduction to $\F$, $\omega_1^{f(M)} \neq \alpha$ (since one assumed that $\omega_1^{f(P)} = \alpha$, if this could occur among the three models). Thus, $\omega_1^{f(M)} < \alpha$ since $\beta$ is the next admissible ordinal after $\alpha$. Observe that
$$X \in [M]_{E_{\cong_\SCRL}^\varphi} \Leftrightarrow f(X) \in [f(M)]_\F$$
Let $y \in \cantorspace$ be such that $z \leq_T y$ and $\omega_1^y = \alpha$ (which exists due to Theorem \ref{Sacks theorem}). $[f(M)]_{\F}$ is $\lborel(y)$ by Proposition \ref{F classes borel in real of higher church kleene}. This shows that $[M]_{\cong_\SCRL}$ is $\Sigma_1^1(y,z) = \Sigma_1^1(y)$. $\omega_1^y = \alpha < R(M)$. This contradicts Theorem \ref{isomorphism class not analytic in some real of small omega1}.

So it has been shown that $\omega_1^{f(M)} > \alpha$. But since $\beta$ is the smallest admissible ordinal greater than $\alpha$, $\omega_1^{f(M)} \geq \beta$. The same exact argument shows $\omega_1^{f(N)} \geq \beta$. 

By Fact \ref{model of high omega1}, let $\CA$ be a countable admissible set extending $L_\beta(z)$ containing $\text{CSS}(M)$ and $\text{CSS}(N)$ with $o(\CA) = \beta$. 

Since $f$ is $\lborel(z)$, let $U$ be a $z$-recursive tree on $2 \times 2 \times \omega$ such that for all $X \in S(\SCRL)$ and $r \in \cantorspace$,
$$(X,r) \in f \Leftrightarrow [U^{(X,r)}] \neq \emptyset$$

Let $\SCRU$ be a language consisting of:

(i) A binary relation symbol $\dot\in$.

(ii) For each $e \in A$, a constant symbol $\bar{e}$.

(iii) Six distinct symbols $\dot R$, $\dot S$, $\dot c$, $\dot d$, $\dot u$, and $\dot v$. 

\noindent $\SCRU$ may be considered as a $\Delta_1$ definable language in $A$. 

Let $T$ be the theory in the countable admissible fragment $\SCRU_{\CA}$ consisting of the following sentences:

(I) $\KP$

(II) For each $e \in A$, $(\forall v)(v \dot\in \bar{e} \Rightarrow \bigvee_{z \in e} v = \bar{z})$. 

(III) $\dot R, \dot S, \dot c, \dot d \subseteq \bar{\omega}$. $\dot u$ and $\dot v$ are functions from $\bar{\omega} \rightarrow \bar{\omega}$. 

(IV) $\dot u \in [U^{(\dot R, \dot c)}]$ and $\dot v \in [U^{(\dot S, \dot d)}]$. 

(V) $\dot R \models \text{CSS}(M)$ and $\dot S \models \text{CSS}(N)$. 

(VI) For all $\xi < \beta$, $\bar{\xi}$ is not admissible in $\dot c$ and $\bar{\xi}$ is not admissible in $\dot d$. 

\noindent $T$ may be considered a $\Sigma_1$ definable theory in $\CA$. Note that $\CA$ was chosen so that (V) would be expressible.

Since $(M, f(M)) \in f$ and $(N, f(N)) \in f$, let $u, v \in \bairespace$ be such that $u \in [U^{(M,f(M))}]$ and $v \in [U^{(N, f(N))}]$. 

Now to show $T$ is consistent: Consider the following $\SCRU$-structure $\mathcal{G}$: The domain of $\mathcal{G}$ is $G = H_{\aleph_1}$. For each $e \in A$, $\bar{e}^\mathcal{G}= e$. $\dot R^\mathcal{G}= M$. $\dot S^\mathcal{G}= N$, $\dot c^\mathcal{G} = f(M)$, $\dot d^\mathcal{G} = f(N)$, $\dot u^\mathcal{G} = u$, and $\dot v^\mathcal{G} = v$. Then $\mathcal{G} \models T$. 

By Theorem \ref{solid model existence theorem}, there exists a solid model $\mathcal{H} \models T$ with $o(\mathcal{H}) = o(\CA)$. Let $R = \dot R^\mathcal{H}$ and $S = \dot S^\mathcal{H}$. Then $f(R) = \dot c^\mathcal{H}$ and $f(S) = \dot d^\mathcal{H}$ since $\dot u^\mathcal{H} \in [U^{(R, \dot c^\mathcal{H})}]$ and $\dot v^\mathcal{H} \in [U^{(S, \dot d^\mathcal{H})}]$. As in the proof of Sacks' theorem, $\omega_1^{f(R)} = \omega_1^{f(S)} = \beta$. By the absoluteness of satisfaction (from $\mathcal{H}$ to $\text{WF}(\mathcal{H})$ to $V$), $R \models \text{CSS}(M)$ and $S \models \text{CSS}(N)$. Hence in $V$, $R$ and $S$ are not $\SCRL$-isomorphic. However, $\omega_1^{f(R)} = \omega_1^{f(S)} = \beta$ implies that $f(R) \ \F \ f(S)$. This contradicts $f$ being a reduction. 

This proves the theorem for $\alpha \in \Lambda(y) \cap \omega_1$. The statement $f$ witnesses $E_\SCRL^\varphi \leq_\borel \F$ is $\bPi_2^1$. So the same argument as at the end of the proof of Theorem \ref{spectrum of aleph0 borel reduction Ew1 to F} shows the results holds for all $\alpha \in \Lambda(y)$. 
\end{proof}

\Begin{corollary}{constructible universe no reduction counterexample to f}
In $L$ (and any set generic extension of $L$), there is no recursive language $\SCRL$ and counterexample $\varphi \in \SCRL_{\omega_1\omega}$ such that $E_\SCRL^\varphi \leq_\borel \F$. 
\end{corollary}

\begin{proof}
There is no $z \in \cantorspace$ having the property of Theorem \ref{counterexample vaught not reducible F} in $L$ or set generic extensions of $L$. 
\end{proof}

\Begin{corollary}{sy friedman question answer in L}
In $L$ (and set generic extensions of $L$), Question \ref{sy friedman vaught question} has a negative answer. 
\end{corollary}

\begin{proof}
This follows from Corollary \ref{sy friedman question answer in L} and the remarks following Question \ref{sy friedman vaught question}.
\end{proof}

This leaves open whether there is an answer to Question \ref{sy friedman vaught question} in $\mathsf{ZFC}$.

\bibliographystyle{amsplain}
\bibliography{references}

\end{document}